\documentclass{article}
\usepackage{amsmath,amsfonts}
\usepackage[francais]{babel}

\makeatletter

\newenvironment{preuve}[1][\hskip-1ex]{\par\reset@font{\slshape %
    Preuve\hskip1ex#1\pointir}}{\hfill\finpr\par}
\newenvironment{remarques}{\par %
\noindent\textbf{Remarques}\pointir}{\hfill\par}
\newenvironment{remarque}{\par %
\noindent\textbf{Remarque}\pointir}{\hfill\par}  
\newenvironment{exemple}{\par %
\noindent\textbf{Exemple}}{\hfill\par}

\makeatother

\newtheorem{theo}{theoreme}[section]
\newtheorem{theoreme}[theo]{Th\'eor\`eme}
\newtheorem{proposition}[theo]{Proposition}
\newtheorem{lemme}[theo]{Lemme}
\newtheorem{corollaire}[theo]{Corollaire}
\newtheorem{definition}[theo]{Definition}
\newtheorem{question}[theo]{Question}

\long\def\InsertFig#1 #2 #3 #4\EndFig{\hbox{\hskip #1 mm$
\vbox to #2mm{\vfil\includegraphics{#3}}#4$}}
\long\def
\LabelTeX#1 #2 #3\ELTX{\rlap{\kern#1mm\raise#2mm\hbox{#3}}}%

\def\build#1#2\fin{\mathrel{\mathop{\kern0pt#1}\limits#2}}
\def\pointir{.\kern.4em\ignorespaces}%
\def\openbox{\leavevmode
  \hbox to.77778em{%
  \hfil\vrule
  \vbox to.675em{\hrule width.6em\vfil\hrule}%
  \vrule\hfil}}
\def\finpr{\openbox}

\def\hfl{{\hbox to 8mm{\rightarrowfill}}}

\def\hf#1{\hhf^{n_{#1}}}
\def\sn#1{S^{n_{#1}{-}1}}
\newcommand{\xgg}{(\wt X,g^1_0\oplus g^2_0)}

\def\Ent{\mathop{\rm Ent}\nolimits}
\def\ci{{\cal C}^\infty}
\def\limti{\dis\lim_{t\to +\infty}}
\def\pa{\partial}

\def\Ker{\mathop{\rm Ker}\nolimits}
\def\ric{\mathop{\rm Ricci}\nolimits}

\def\vol{\mathop{\rm vol}\nolimits}
\def\det{\mathop{\hbox{\rm d\'et}}\nolimits}
\def\deg{\mathop{\hbox{\rm deg}}\nolimits}

\def\minvol{\mathop{\rm minvol}\nolimits}

\def\divv{\mathop{\rm div}\nolimits}

\def\Jac{\mathop{\rm Jac}\nolimits}

\def\Isom{\mathop{\rm Isom}\nolimits}
\def\barr{\mathop{\rm bar}\nolimits}
\def\tr{\mathop{\rm trace}\nolimits}

\def\PO{\mathop{\rm PO}\nolimits}

\def\bs{\overline{s}}
\def\bg{\overline{g}}
\def\bF{\overline{F}}
\def\btheta{\overline{\theta}}

\def\ol{\overline}

\newcommand\zzf{\mathbf{Z}}
\providecommand{\cad}{c'est-\`a-dire }%
\newcommand\bc{\mathcal{B}}
\newcommand\vc{\mathcal{V}}

\newcommand\mc{\mathcal{M}}
\newcommand{\ld}{, \ldots,}
\newcommand{\prodd}{\mathop{\prod}\limits}
\newcommand{\summ}{ \mathop{\sum}\limits}
\newcommand{\opp}{ \mathop{\bigoplus}\limits}
\newcommand{\cupp}{ \mathop{\bigcup}\limits}

\def\ixg{\Isom(\wt X, \tilde g_0)}
\providecommand{\cf}{{\upshape cf. }}%

\providecommand{\resp}{{resp. }}%
\newcommand{\wt}{\widetilde}
\newcommand{\wh}{\widehat}

\providecommand{\la}{\longrightarrow}
\newcommand\hhf{\mathbf{H}}
\newcommand\rrf{\mathbf{R}}
\newcommand\hg{\mathfrak{h}}
\newcommand\ggg{\mathfrak{g}}
\newcommand{\dis}{\displaystyle}

\begin{document}

\title{In\'egalit\'es de Milnor-Wood g\'eom\'etriques}
\author{G. Besson, G. Courtois et  S. Gallot}

\maketitle

\section{Introduction}
La c\'el\`ebre in\'egalit\'e de Milnor-Wood (\cite{Mil} et \cite{Woo})
affirme que, si 
$$E\longrightarrow \Sigma$$
est un fibr\'e plat en fibres $S^1$ sur la surface compacte $\Sigma$ de genre
$\gamma\geq 2$, alors la caract\'eristique d'Euler de ce fibr\'e, not\'ee
$\chi (E)$ v\'erifie,
$$|\chi (E)|\leq |\chi (\Sigma )|=2\gamma -2\, ,$$
l'\'egalit\'e ayant lieu si $E$ est le fibr\'e tangent de $\Sigma$.

Un fibr\'e plat \'etant d\'efini par une repr\'esentation de $\Pi_1(\Sigma )$,
 l'in\'egalit\'e ci-dessus est en fait une restriction impos\'ee \` a 
cette repr\'esentation \`a valeurs dans le groupe des hom\'eomorphismes du 
cercle. Les valeurs possibles de la caract\'eristique d'Euler de $E$ sont 
d\'ecrites dans \cite{Gol}. 

Dans cette article nous envisageons une g\'en\'eralisation, 
en dimension sup\'erieure de cette in\'egalit\'e. Pour cela nous d\'efinissons
le volume d'une repr\'esentation. Plus pr\'ecis\'ement, soit $M$ une 
vari\'et\'e diff\'erentielle ferm\'ee de dimension $n$ et soit $\rho$ une 
repr\'esentation de son groupe fondamental dans le groupe d'isom\'etrie 
d'une vari\'et\'e sym\'etrique de courbure n\'egative de dimension $n$ 
et simplement connexe, not\'ee $\wt X$. Consid\'erons une application 
$\tilde{f}$ du rev\^etement universel de $M$ dans $\wt X$ \'equivariante
par rapport \`a $\rho$, alors, si $\omega$ d\'esigne la forme volume de
$\wt X$, la forme $\tilde{f}^*\omega$ passe au quotient sur $M$.
\begin{definition}
On appelle volume de la repr\'esentation $\rho$ le nombre,
$$\vol (\rho )=\int_M\tilde{f}^*\omega\,.$$
\end{definition}
Dans certains cas ce nombre peut-\^etre interpr\'et\'e comme la classe d'Euler
d'un fibr\'e plat. Des bornes sup\'erieures de $\vol (\rho )$ existent. Elles
reposent souvent sur le choix d'une famille de sections particuli\`eres du
fibr\'e plat. Dans \cite{Cor}, par exemple, K.~Corlette utilise des sections 
harmoniques pour d\'emontrer un th\'eor\` eme de rigidit\'e sur les 
repr\'esentations de volume maximal. Le cas o\`u $\wt X$ est hyperbolique r\'eel est
abord\'e par A.~Reznikov dans \cite{Rez1}; l'auteur y prouve une in\'egalit\'e
optimale et c'est ce type de r\'esultats que nous \'etendons dans le 
pr\'esent travail. Le cas d'\'egalit\'e dans l'in\'egalit\'e de A.~Reznikov
est prouv\'e par N.~Dunfield dans \cite{Dun} et dans \cite{Cou}, il consiste \`a montrer que,
si le volume est maximal, la repr\'esentation est fid\`ele et discr\`ete. Signalons l'article \cite{Klin} dans lequel l'auteur d\'ecrit une autre notion de volume de repr\'esentations et construit de nouveaux invariants num\'eriques.

 Dans le cas o\`u $\wt X$ est l'espace hyperbolique r\'eel nous prouvons, dans cet article,
que le volume des repr\'esentations est constant sur les composantes
connexes de l'espace des repr\'esentations. C'est un r\'esultat
\'evident
 lorsque la dimension est paire car, dans ce cas, le volume
est aussi un
 nombre d'Euler, mais nouveau dans le cas de dimension
impaire. Plus pr\'ecis\'ement nous prouvons le

\begin{theoreme}\label{vol-constant}
Soit $M$ une vari{\'e}t{\'e} diff{\'e}rentielle ferm{\'e}e et orient{\'e}e et 
$\rho_t : \Pi_1(M)\longrightarrow Isom({\widetilde X})$ une famille de repr{\'e}sentations
qui d{\'e}pend de mani{\`e}re $C^1$ du param{\`e}tre $t\in \bf{R}$, alors le volume 
$\vol (\rho _t)$ est constant.
\end{theoreme}
La m\'ethode employ\'ee
consiste \` a utiliser la formule de Schl\"afli (voir aussi \cite{Bon2}). Il
s'agit d'une approche nouvelle dans ce contexte; en fait nous construisons 
un ``poly\`edre'' g\'eod\'esique hyperbolique dans $\widetilde X$ 
\'equivariant par rapport
\`a l'image de $\rho$. Il s'agit d'une r\'eunion de simplexes hyperboliques
g\'eod\'esiques invariants par l'image de $\rho$; les simplexes pouvant 
se chevaucher ils ne fournissent   pas une triangulation de $\widetilde X$. 
Nous construisons
ensuite une application $\rho$-\'equivariante poly\`edrale. Ceci permet
alors de calculer le volume de la repr\'esentation. La formule de
Schl\"afli ainsi qu'un peu de th\'eorie du degr\'e permet
alors de montrer la constance de ce volume.
Un corollaire frappant est une preuve tr\`es simple du r\'esultat suivant 
de T.~Soma (\cite{Som}),
\begin{theoreme}[T.~Soma]
Soit $Y$ une vari\'et\'e diff\'erentielle ferm\'ee de dimension $3$.
L'ensemble des vari\'et\'es hyperboliques ferm\'ees $X$, de dimension $3$ 
telles qu'il existe une application continue de degr\'e non nul de $Y$ sur
$X$, est fini.
\end{theoreme}
La preuve se r\'esume comme suit. Appelons $f$ l'application de degr\'e non nul
de $Y$ sur $X$ et $\rho$ la repr\'esentation induite de $\Pi_1(Y)$ dans $\Pi_1(X)$, 
alors $\vol (\rho )=\deg (f)\vol (X)$. Le th\'eor\`eme \ref{principal-intro}
affirme qu'il existe une constante $C(M)$ telle que $\vol (\rho)\leq C(M)$. Par ailleurs 
le volume des vari\'et\'es hyperboliques ferm\'ees est minor\'e par une constante 
universel (qui d\'epend de la dimension). Ceci montre que le degr\'e de $f$ ne peut
prendre qu'un nombre fini de valeurs. Le volume de la repr\'esentation $\rho$
\'etant constant sur les composantes connexes de l'espace des repr\'esentations
de $\Pi_1(Y)$ dans $\Pi_1(X)$ il ne prend \'egalement qu'un nombre fini de valeurs.
Finalement le volume de la vari\'et\'e hyperbolique $X$ ne prend qu'un nombre fini de valeurs.
Un r\'esultat de W.~Thurston affirme alors qu'il n'y a qu'un nombre fini de vari\'et\'es
$X$ possibles. Le lecteur peut consulter les d\'etails dans le paragraphe \ref{applications}.

Le th\'eor\`eme \ref{vol-constant} est en fait un corollaire d'un r\'esultat plus g\'en\'eral, que
nous d\'ecrivons maintenant.
Rappelons la d\'efinition de l'entropie volumique d'une vari\'et\'e 
Riemannienne $(Y,g)$. Pour $x\in \wt Y$, on d\'esigne par $B(x,R)$ la 
boule g\'eod\'esique de centre $x$ et de rayon $R$, alors on d\'efinit
$$\Ent(Y,g)=\lim_{R\to \infty}\frac{1}{R}\log (\vol (B(x,R)))\,.$$
Dans ce qui suit $\wt X=\prod_1^p\wt X_i$ est le produit des espaces sym\'etriques simplement
connexes de courbure strictement n\'egative, $\wt X_i$. Chacune des vari\'et\'es $\wt X_i$ est munie 
d'une m\'etrique $\alpha_ig_0^i$, o\`u $g_0^i$ est sym\'etrique normalis\'ee (de courbure comprise entre $-4$ et $-1$, par exemple) et $\alpha_i$ est un r\'eel strictement positif. Parmi tous les choix de nombres $\alpha_i$ il en est un qui donne une entropie volumique minimale (voir la proposition \ref{min}); nous noterons $g_0$ la m\'etrique correspondante sur $\wt X$ qui est de dimension $n$. Son entropie volumique est un 
nombre calculable. Nous prouvons,
\begin{theoreme}\label{principal-intro}
Soit $Y$ une vari\'et\'e riemannienne ferm\'ee de dimension $n$ et $\rho$
une repr\'esentation de $\Pi_1(Y)$ dans $\Isom (\wt X)$, alors
\begin{itemize}
\item[i)] $\vol (\rho )\leq \big (\frac{\Ent (Y,g)}{\Ent (\wt X, g_0)}\big )^n
\vol(Y,g)\,.$
\item[ii)] L'\'egalit\'e dans l'in\'egalit\'e ci-dessus a lieu si, et
seulement si, la repr\'esentation $\rho$ est injective, $X=\wt X/\rho (\Pi_1
(Y))$ est une vari\'et\'e compacte et $(Y,g)$ est homoth\'etique \`a 
$(X,g_0)$
\end{itemize}
\end{theoreme}  
Ce r\'esultat \'etait annonc\'e en 1997 dans \cite{Gui} et \'enonc\'e en 
1998 dans \cite{Cou}.
Il g\'en\'eralise le cas o\`u la repr\'esentation a une image discr\`ete et
cocompacte, c'est-\`a-dire l'analogue des th\'eor\`emes de \cite{BCG1} pour
le cas o\`u l'espace localement sym\'etrique compacte est localement un
produit d'espaces sym\'etriques de rang 1. Ce dernier r\'esultat, concernant 
les produits d'espaces sym\'etriques de rang 1 avec image discr\`ete 
cocompacte, est \'enonc\'e par Ch.~Connell et B.~Farb dans \cite{Con-Far}.

La preuve de l'in\'egalit\'e se fait en exhibant une famille
d'applications 
 $\rho$-\'equivariantes de $\wt Y$ sur $\wt X$
construites par la m\'ethode
 introduite dans \cite{BCG1}.  Le cas
d'\'egalit\'e est beaucoup plus difficile car l'image de $\rho$ n'est
pas suppos\'ee discr\`ete; plus pr\'ecis\'ement, nous montrons que,
dans le cas d'\'egalit\'e, la famille d'applications
$\rho$-\'equivariantes que nous construisons converge vers une
application harmonique; ceci permet, en particulier, de montrer que la
limite est de classe $C^\infty$. La combinaison des propri\'et\'es
li\'ees \`a l'harmonicit\'e et de celles li\'ees \`a la construction
ci-dessus conduit au r\'esultat.

Remarquons que les applications $\rho$-\'equivariantes construites
sont particuli\`erement adapt\'ees \`a l'\'etude du volume et
conduisent \`a des r\'esultats optimaux comparables, dans un cadre plus
g\'en\'eral, \`a ceux de N.~Dunfield \cite{Dun}. Signalons \'egalement
un travail r\'ecent de S.~Francaviglia et B.~Klaff \cite{Fra-Kla} dans
lequel les auteurs utilisent une int\'eressante variante de la
construction de \cite{BCG2} pour \'etudier le cas o\`u $Y$ est une
vari\'et\'e hyperbolique de volume fini.  
 
Enfin, l'in\'egalit\'e
ci-dessus peut s'interpr\'eter agr\'eablement dans le
 cadre de la
cohomologie born\'ee (voir \cite{Gro}). Le r\'ecent travail de
M.~Burger, A.~Iozzi et A.~Wienhard (\cite{Bur-Ioz-Wie}) d\'eveloppe ce
point 
 de vue et aboutit \`a de tr\`es jolis r\'esultats concernant
les 
 repr\'esentations du groupe fondamental des surfaces.

Nous tenons \`a remercier A.~Reznikov, M.~Boileau, D.~Cooper et S.~Francaviglia pour leur aide et leurs commentaires lors de la redaction de cet article.

\section{G{\'e}om{\'e}trie des espaces produits}

{\`A} titre d'exemple, nous d{\'e}crirons la g{\'e}om{\'e}trie de l'espace $\xgg= (\hf
1\times \hf 2,g^1_0 \oplus g^2_0)$ muni 
de la m{\'e}trique produit o{\`u}  $(\hhf^{n_1},g^1_0)$ ($\resp
(\hf 2,g^2_0)$) d{\'e}signe l'espace hyperbolique simplement connexe
de dimension $n_1$ ($\resp n_2$) (de courbure constante {\'e}gale {\`a}
$-1$). Pour un expos{\'e} g{\'e}n{\'e}ral sur les espaces sym{\'e}triques, nous
renvoyons {\`a} \cite{Hel}.

\subsection{G{\'e}od{\'e}siques}

Soient $x = (x_1,x_2) \in \wt X$ et $ u = (u_1,u_2) \in T_{(x_1,x_2)}
\wt X$ tels que $\|u\|^2_{g^1_0 \oplus g^2_0} = \|u_1\|^2_{g^1_0} +
\|u_2\|^2_{g^2_0} = 1$, alors la g{\'e}od{\'e}sique de $X$, not{\'e}e $c_u$, partant
de $x$ et de vitesse initiale $u$ est $c_u(t) = (c_1(t),c_2(t))$, o{\`u}
 $c_i$ ($i=1,2$) est la g{\'e}od{\'e}sique de $\hf i$ partant de $x_i$ et de
vitesse initiale $u_i$. Une g{\'e}od{\'e}sique d{\'e}finie par un vecteur $u =
(u_1,u_2)$ telle que $u_1 =0$ ou bien $u_2=0$ est dite singuli{\`e}re~;
ces cas correspondent {\`a}
$$c_u(t) = (x_1,c_2(t)) \hbox{~~ou~~} c_u(t) = (c_1(t),x_2)~.$$
Une g{\'e}od{\'e}sique d{\'e}finie par un vecteur $u = (u_1,u_2)$ tel que $u_i \ne
0$, pour $i = 1,2$, est dite r{\'e}guli{\`e}re.

\subsection{Courbures et plats}

La courbure sectionnelle de $\xgg$, qui se calcule
ais{\'e}ment, est n{\'e}gative ou nulle. Soit alors $x = (x_1,x_2) \in X$, $u
= (u_1,u_2) \in T_x\wt X$, un vecteur r{\'e}gulier, alors l'application
\begin{align*}
\rrf^2 &\la \wt X\\
(t,s) & \longmapsto \big(c_1(t/\alpha_1 ),c_2(s/\alpha_2 )\big)
\end{align*}
o{\`u} $\alpha_1  = \|u_1\|_{g^1_0}$ et $\alpha_2  = \|u_2\|_{g^2_0}$ r{\'e}alisent un
plongement isom{\'e}trique de $\rrf^2$ muni de sa m{\'e}trique euclidienne
dans $\xgg$. On peut v{\'e}rifier par le calcul que l'image de cette
application est totalement g{\'e}od{\'e}sique (voir \cite{Hel}, pp.~~~) ou bien
constater que, si $\sigma _i$ d{\'e}signe la sym{\'e}trie orthogonale par
rapport {\`a} la g{\'e}od{\'e}sique $c_i$  dans $(\hf i, g^i_0)$, l'image de
l'application ci-dessus est l'ensemble des points fixes de $\sigma _1
\times \sigma _2$ dans $\wt X$~; il s'agit donc d'un sous-espace
totalement g{\'e}od{\'e}sique plat et qui est, de plus, de dimension maximale
avec ces propri{\'e}t{\'e}s~: $\xgg$ est un espace sym{\'e}trique de rang 2. 
{\bf Nous noterons d{\'e}sormais $\bar g_0$ la m{\'e}trique $g^1_0 \oplus g^2_0$}.

\begin{remarque}D'une mani{\`e}re g{\'e}n{\'e}rale, si $\wt X$ est le produit
riemannien de $p$ espaces sym{\'e}triques de courbure strictement
n{\'e}gative, alors $X$ est de rang $p$.
\end{remarque}

\subsection{M{\'e}triques localement sym{\'e}triques}

On peut munir la vari{\'e}t{\'e} diff{\'e}rentielle $\wt X$ d'autres m{\'e}triques
localement sym{\'e}triques~; en effet, pour $\alpha_1 $ et $\alpha _2 $ deux
nombres r{\'e}els strictement positifs, on d{\'e}finit~:
$$g_{\alpha_1 ,\alpha _2 } = \alpha_1 ^2 g^1_0 \oplus \alpha _2 ^2g^2_0~.$$
Contrairement aux espaces sym{\'e}triques irr{\'e}ductibles, les espaces
sym{\'e}triques produits sont flexibles.

\subsection{Groupe d'isom{\'e}tries}

On d{\'e}termine ais{\'e}ment le groupe d'isom{\'e}tries de $(\wt X,g_{\alpha_1
,\alpha _2 })$. En effet, si $n_1 \ne n_2$
$$\Isom (\wt X,g_{\alpha_1,\alpha _2 }) = \Isom (\hf 1,g^1_{0})\times
\Isom (\hf 2,g^2_{0})~.$$ 
Si $n_1=n_2$ et $\alpha_1=\alpha_2$, l'{\'e}change des deux facteurs est une isom{\'e}trie
suppl{\'e}mentaire qui est involutive~; le groupe d'isom{\'e}trie de 
$(\wt X,g_{\alpha_1,\alpha_1 })$ est donc une extension  de $\zzf/2\zzf$
par le groupe 
$\Isom (\hf 1,g_{0})\times \Isom (\hf 2,g_{0})$.

\subsection{Fonctions de Busemann}

On rappelle que, si $(M,g)$ est une vari{\'e}t{\'e} riemannienne compl{\`e}te et si
$c : \rrf \to M$ est une g{\'e}od{\'e}sique minimisante sur toute sa longueur et
param{\'e}tr{\'e}e par l'abscisse curviligne (c'est-{\`a}-dire, $c$ est un plongement
isom{\'e}trique), alors on d{\'e}finit la fonction de Busemann associ{\'e}e {\`a} $c$,
$$B_c(x) = \limti d(x,c(t)) - t = \limti
\big(d(x,c(t)) - d(c(0),c(t))\big)~.$$
On montre que la limite existe (voir \cite{BGS}, p.~23). Si $(M,g)$ est une
vari{\'e}t{\'e} simplement connexe de courbure n{\'e}gative ou nulle son bord {\`a}
l'infini (voir \cite{BGS}, p.~) s'identifie {\`a} une sph{\`e}re de dimension
$n-1$, o{\`u} $n =\dim M$, gr{\^a}ce au choix d'un point $O\in M$ qui sert
d'origine. Chaque point $\theta  \in \pa M$, le bord {\`a} l'infini de
$M$, d{\'e}termine une g{\'e}od{\'e}sique minimisante sur toute sa longueur, {\`a}
savoir, l'unique g{\'e}od{\'e}sique $c$ qui passe par $O$ et telle que $\limti
c(t) = \theta $. La fonction de Busemann correspondante est not{\'e}e
$B(\cdot,\theta )$. Remarquons qu'elle d{\'e}pend du choix de l'origine.

Dans notre situation, il est souhaitable de travailler sur une partie
du bord qui refl{\`e}te mieux la structure produit. Pour la vari{\'e}t{\'e} $\wt
X$ ci-dessus le bord {\`a} l'infini s'identifie {\`a} $S^{n_1+n_2-1}$ (pour
toutes les m{\'e}triques $g_{\alpha_1 ,\alpha _2 }$) apr{\`e}s le choix d'une
origine. Nous utiliserons $\sn 1 \times  \sn 2 \subset 
S^{n_1+n_2-1}$ qui  s'identifie dans $\pa \wt X$ {\`a} $\pa \hf 1
\times \pa \hf 2$. Plus pr{\'e}cis{\'e}ment, consid{\'e}rons, par exemple, la
m{\'e}trique $\bar g_{0} = g^1_0 \oplus g^2_0$, appelons $O = (O_1,O_2)$ une
origine de $\wt X = \hf 1 \times \hf 2$, le bord de $\wt X$
s'identifie aux rayons g{\'e}od{\'e}siques param{\'e}tr{\'e}s par longueur d'arc et
partant de $O$~; nous ne consid{\'e}rerons que les g{\'e}od{\'e}siques $c =
(c_1,c_2)$ o{\`u} $c_i$ est une g{\'e}od{\'e}sique de $\hf i$, telle que, pour tout
$t \in \rrf$, $\| \dot c_1(t)\|_{g^1_0} = \|\dot c_2(t)\|_{g^2_0}$~; nous les
appellerons g{\'e}od{\'e}siques diagonales. Elles sont donc param{\'e}tr{\'e}es par un
point $\theta  = (\theta _1,\theta _2)$ o{\`u} $\theta _i \in \sn i =
\pa \hf i$. Il s'agit du bord de Furstenberg (voir ~~~), mais nous
n'utiliserons pas sa description probabiliste. Nous le noterons $\pa_F
\wt X$. Il est important de noter que nous utiliserons toujours ce bord; en effet, si nous changeons la m{\'e}trique en $g_{\alpha_1 ,\alpha_2 }$, nous pouvons consid{\'e}rer des $g_{\alpha_1 ,\alpha_2
}$-g{\'e}od{\'e}siques $c = (c_1,c_2)$ telles que ${1\over \alpha_1 } \|\dot
c_1(t)\|_{g_{\alpha_1} } = {1\over \alpha_2 } \|\dot c_2
(t)\|_{g_{\alpha_2 }}$, 
o{\`u} $g_{\alpha_i}  = \alpha_i ^2 g^i_0$; elles d\'efinissent un bord qui s'identifie \`a $\pa_F \wt X$. 

\begin{remarque}
Lorsque $n_1 = n_2 =2$ et $\alpha_1  = \alpha_2  = 1$, le bord de
Furstenberg de $\wt X$, $S^1 \times S^1 \subset S^3 = \pa \wt X$,
s'identifie naturellement {\`a} un tore de Clifford dans $S^3$.
\end{remarque}

Maintenant, pour $\theta =(\theta
_1,\theta _2) \in \sn 1 \times \sn 2$, on note $\overline B_0(\cdot,\theta )$ la fonction de Busemann de
$(X,\bar g_{0})$ correspondante (l'origine $O=(O_1,O_2)$ {\'e}tant fix{\'e}e), et
$B_i (\cdot, \theta _i)=B_{O_i} (\cdot, \theta _i)$, $i=1,2$, la fonction de Busemann de $(\hf
i,g^i_0)$, on a~:

\begin{lemme}\label{1.1}Avec les notations ci-dessus, si $x = (x_1,x_2)\in \wt
X$
$$\overline B_0(x,\theta ) = {1\over \sqrt{2}} \big(B_1(x_1,\theta
_1)+B_2(x_2,\theta _2)\big)~.$$
\end{lemme}

\begin{preuve} Soit $c$ la g{\'e}od{\'e}sique param{\'e}tr{\'e}e par l'abscisse curviligne
d{\'e}finie par $\theta $ et telle que $c(0) = O = (O_1,O_2)$. Alors, si $c
= (c_1,c_2)$, on a 
$\|\dot c_1\| = {1\over \sqrt 2} = \|\dot c_2\|$, d'o{\`u} 
$$d_i(x_i,c_i(t)) = {1\over \sqrt 2} t + B_i(x_i,\theta _i) +
\varepsilon _i(t),~~ i=1,2$$
avec $\varepsilon _i (t) \build \hfl_{t\to +\infty}\fin 0$. Ici, $d_i$
d{\'e}signe la distance dans le facteur $i=1,2$.

Le lemme se d{\'e}duit alors du d{\'e}veloppement limit{\'e} de
$$d(x,c(t)) -t = \big( d^2_1(x_1,c_1(t)) + d^2_2
(x_2,c_2(t))\big)^{1/2} - t~.$$
\end{preuve}

De m{\^e}me, si $B_{\alpha_1 ,\alpha_2 }(\cdot, \theta )$ d{\'e}signe la fonction de
Busemann de $(\wt X,g_{\alpha_1 ,\alpha_2 })$ o{\`u} $\theta $ est dans le bord
d{\'e}fini ci-dessus, on a~:

\begin{lemme} Avec les notations ci-dessus, si $x = (x_1,x_2) \in \wt
X$ 
$$B_{\alpha_1 ,\alpha_2 }(x,\theta ) = {1\over \sqrt {\alpha_1 ^2+\alpha_2 ^2}}
\big( \alpha_1 B_1(x_1,\theta _1) + \alpha_2 B_2(x_2,\theta
_2)\big)~.$$
\end{lemme}

La preuve de ce lemme se fait comme celle du lemme \ref{1.1}.

\subsection{{\'E}l{\'e}ment de volume}

Si on note $dv_g$ l'{\'e}l{\'e}ment de volume d'une m{\'e}trique riemannienne $g$,
il est imm{\'e}diat que
$$dv_{g_{\alpha_1 ,\alpha_2 }} = \alpha_1 ^{n_1} \alpha_2 ^{n_2}
dv_{g^1_0}\otimes dv_{g^2_0}$$
o{\`u} $dv_{g^i_0}$ d{\'e}signe l'{\'e}l{\'e}ment de volume de $(\hf i,g^i_0)$ pour $i =
1,2$.

\subsection{Entropie}

On rappelle la d{\'e}finition de l'entropie (volumique) d'une vari{\'e}t{\'e}
riemannienne $(M,g)$ que nous supposerons compacte pour
simplifier. Soit $m \in \wt M$ un point du rev{\^e}tement universel $\wt
M$ de $M$ alors la quantit{\'e} suivante existe et ne d{\'e}pend pas de $x$, 
$$\Ent(g) = \lim_{R\to +\infty} {1\over R} \log \big(\vol
(B_{\wt M}(x,R))\big)$$
o{\`u} $B_{\wt M} (x,R)$ d{\'e}signe la boule m{\'e}trique de centre $x$ et de
rayon $R$ dans $\wt M$ muni de la m{\'e}trique relev{\'e} de $g$.

Par d{\'e}finition $\Ent(g)$ est l'entropie de la vari{\'e}t{\'e} riemannienne
$(M,g)$, elle ne d{\'e}pend de $M$ qu'{\`a} travers la relev{\'e}e de $g$ {\`a} $\wt
M$. Par abus de langage, nous parlerons de l'entropie de $g_{\alpha_1
,\alpha_2 }$ sur $\wt X$.

\begin{proposition}
Pour tous $\alpha_1 ,\alpha_2 $ positifs
$$\Ent(g_{\alpha_1 ,\alpha_2 }) = \sqrt{{(n_1{-}1)^2\over \alpha_1 ^2} +
{(n_2{-}1)^2\over \alpha_2  ^2}}~.$$
\end{proposition}

\begin{preuve} Le calcul de l'entropie des espaces sym{\'e}triques est fait dans
\cite{BCG1}. Rappelons que l'entropie d'un produit v{\'e}rifie
$$\Ent(g_{\alpha_1 ,\alpha_2 })^2 = {\Ent (\hf 1,g^1_0)^2 \over \alpha_1 ^2} +
{\Ent(\hf 2,g^2_0)^2 \over \alpha_2 ^2}~.$$
\end{preuve}

Dans cet article on se propose de prouver un th{\'e}or{\`e}me d'entropie
minimale (voir l'introduction) \cad de minimum de l'entropie {\`a} volume
fix{\'e}. Dans ce paragraphe nous examinons cette question pour la
famille de m{\'e}trique $g_{\alpha_1 ,\alpha_2 }$. Plus pr{\'e}cis{\'e}ment, soit
$\Gamma $ un sous-groupe discret cocompact de $\Isom(\hf 1,g^1_0) \times
\Isom(\hf 2,g^2_0)$, agissant sans points fixes sur $\wt X$. Ce groupe
agit par isom{\'e}tries sur $\wt X$ pour toutes les m{\'e}triques $g_{\alpha_1
,\alpha_2 }$, on peut donc munir le quotient $X = \wt X/\Gamma $ des
m{\'e}triques induites que nous noterons encore $g_{\alpha_1 ,\alpha_2 }$. Par
ailleurs,
$$\vol(X,g_{\alpha_1 ,\alpha_2 }) = \alpha_1 ^{n_1} \alpha_2 ^{n_2}
\vol(X,\bar g_{0})~.$$

\begin{proposition}
\label{min}
Pour tous $\alpha_1 ,\alpha_2 $ strictement positifs
tels que
$ \alpha_1 ^{n_1} \alpha_2 ^{n_2}=1$, on~a
\begin{eqnarray*}
\Ent(g_{\alpha_1 ,\alpha_2 }) & \ge & \sqrt{n_1+n_2} \left(
\Big({n_1{-}1\over \sqrt{n_1}}\Big)^{n_1}\ \Big({n_2{-}1\over
\sqrt{n_2}}\Big)^{n_2}\right)^{1\over n_1+n_2}\\
& = & \Ent \big(g_{a_1 ,a_2 }\big) \\
\end{eqnarray*}
o{\`u}
$a_1  = \left[\Big( {(n_1{-}1)\ \sqrt{n_2}\over \sqrt{n_1}\
(n_2{-}1)}\Big)^{n_2} \right]^{1\over n_1+n_2}$
,~~ $a_2  = \left[\Big( {(n_2{-}1)\ \sqrt{n_1}\over \sqrt{n_2}\
(n_1{-}1)}\Big)^{n_1} \right]^{1\over n_1+n_2}$.

L'{\'e}galit{\'e}, dans l'in{\'e}galit{\'e} ci-dessus, a lieu si et
seulement si $\alpha_i  = a_i$. 
\end{proposition}

\begin{remarque}
Lorsque les espaces sym\'etriques sont complexes, quaternioniens ou de Cayley, les calculs sont comparables et sont laiss\'es au lecteur.
\end{remarque}
{Dans la suite nous noterons $g_0$ la
m{\'e}trique $g_{a_1 , a_2}$.

\begin{preuve} 
On a $$\Ent(g_{\alpha_1 ,\alpha_2 })^2 = (n_1+n_2) \left(
{n_1 \Big({n_1-1\over \sqrt{n_1}\alpha_1 }\Big)^2 +n_2 \Big({n_2-1\over
\sqrt{n_2}\alpha_2  }\Big)^2 \over n_1+n_2} \right)$$
la fonction $x \mapsto x^2$ {\'e}tant strictement log-concave
$$
\Ent(g_{\alpha_1 ,\alpha_2 })^2
\ge (n_1+n_2) \Big({n_1{-}1\over \sqrt{n_1}}\Big)^{{2n_1\over
n_1+n_2}}\ \Big({n_2{-}1\over \sqrt{n_2}}\Big)^{{2n_2\over
n_1+n_2}}\ \Big({1\over \alpha_1 ^{n_1}\alpha_2 ^{n_2}}\Big)^{2\over
n_1+n_2}$$
d'o{\`u} le r{\'e}sultat
$$
\Ent(g_{\alpha_1 ,\alpha_2 })
\ge \sqrt{n_1+n_2} \left( \Big({n_1{-}1\over \sqrt{n_1}}\Big)^{n_1}\
\Big({n_2{-}1\over \sqrt{n_2}}\Big)^{n_2}\right)^{1\over 
n_1+n_2}~.$$
De plus, par stricte log-concavit{\'e}, l'{\'e}galit{\'e} n'a lieu que si et
seulement si
$${n_1-1\over \sqrt{n_1} \alpha_1 } = {n_2-1\over \sqrt{n_2}\alpha_2} $$
\cad si $a_i  = \alpha_i$.
\end{preuve}

\begin{remarques}
\begin{itemize}
\item[i)] Si $n_1 = n_2$, alors la m{\'e}trique minimisante est
homoth{\'e}tique {\`a} $\bar g_{0}$ (le facteur d'homoth{\'e}tie
{\'e}tant calcul{\'e} de 
sorte {\`a} avoir un volume 1.

\item[ii)] La courbure de Ricci de la m{\'e}trique $g_{\alpha_1 ,\alpha_2 }$
est
$$\ric(g_{\alpha_1 ,\alpha_2 }) = (n_1-1) g^1_0 \oplus (n_2-1)g^2_0~.$$
La m{\'e}trique $g_{\alpha_1 ,\alpha_2 }$ n'est donc d'Einstein que si
$${n_1-1\over \alpha_1 ^2} = {n_2-1\over \alpha_2 ^2}~.$$
Par cons{\'e}quent, en g{\'e}n{\'e}ral, la m{\'e}trique qui minimise la 
fonctionnelle 
$\Ent$, parmi les 
$g_{\alpha_1 ,\alpha_2 }$, n'est pas d'Einstein. Par contre, elle l'est si et seulement si
$n_1 = n_2$.
\end{itemize}
\end{remarques}

De m{\^e}me, si $X$ est un espace produit g{\'e}n{\'e}ral, c'est-{\`a}-dire, si $(X,\bar g) =
(X_1,g_1) \times \cdots \times (X_p,g_p)$, o{\`u} $(X_k,g_k)$ est un
espace sym{\'e}trique de courbure strictement n{\'e}gative, de dimension $n_k$
et d'entropie not{\'e}e $E_k$, on consid{\`e}re les m{\'e}triques,
$$g_\alpha  = \alpha ^2_1 g_1 \oplus \cdots \oplus \alpha ^2_p g_p$$
o{\`u} $\alpha  = (\alpha _1\ld \alpha _p)$ avec $\alpha _k >0$. Alors, on
a la

\begin{proposition}\label{mingeneral}Pour tous $\alpha _1\ld \alpha _p$ r{\'e}els strictement
positifs tels que $\alpha ^{n_1}_1 \cdots \alpha ^{n_p}_p = 1$, on a
$$\Ent(g_\alpha ) \ge \sqrt n \left(\prodd^p_{i=1} \left({E_i\over
\sqrt{n_i}}\right)^{n_i\over n}\right)$$
o{\`u} $n = n_1 + \cdots + n_p = \dim (X)$.

L'{\'e}galit{\'e}, dans l'in{\'e}galit{\'e} ci-dessus, a lieu si, et seulement si  ,
pour tout $i = 1,2\ld p$
$$\alpha _i = a_i = {E_i \over \sqrt{n_i}} \left(\prodd^p_{k=1}
\left({\sqrt{n_k}\over E_k }\right)^{n_k\over n}\right)~.$$
\end{proposition}

\subsection{Mesure de Patterson-Sullivan}

Sur le rev{\^e}tement universel d'une vari{\'e}t{\'e} de courbure strictement
n{\'e}gative, $(M,g)$, on peut d{\'e}finir une famille de mesures qui est
appel{\'e}e (par abus de langage) la mesure de Patterson-Sullivan. Elle
consiste {\`a} associer {\`a} chaque point $m \in \wt M$ 
(le rev{\^e}tement universel de $M$) une mesure bor{\'e}lienne positive sur
$\pa \wt M$, not{\'e}e $\mu _m$. Cette famille est enti{\`e}rement
caract{\'e}ris{\'e}e par les deux propri{\'e}t{\'e}s suivantes~:

{\it i)\/}~~$\dis{d\mu _m\over d\mu _{m'}}(\theta ) = \exp
\big({-}\Ent(g) (B(m,\theta ) - B(m',\theta ))\big)$ (on a
choisi ici une origine $O \in \wt M$ afin de d{\'e}finir $B$). Cette
propri{\'e}t{\'e} affirme que pour $m\ne m'$ les mesures $\mu _m$ et $\mu
_{m'}$ sont absolument continues l'une par rapport {\`a} l'autre et la
densit{\'e} s'exprime comme ci-dessus.

{\it ii)\/}~~ $\forall \gamma  \in \Isom(\wt M)$, $\gamma $ agit par
hom{\'e}omorphisme sur $\pa \wt M$, et
$$\mu _{\gamma (m)} = \gamma _*(\mu _m)$$
(voir \cite{Kni}).

Dans le cas o{\`u} $\wt M$ est un espace sym{\'e}trique de courbure
n{\'e}gative 
ou nulle (et pas strictement n{\'e}gative) une construction est possible
(voir \cite{Alb}, \cite{Quint} bet \cite{Link}). Dans notre situation, \cad
$$(\wt M,g) = (\wt X,\bar g_{0}) = (\hf 1,g^1_0) \times (\hf 2,g^2_0)$$
la famille de mesures suivante, port\'ees par $\partial H^{n_1}x\partial H^{n_2}$ v{\'e}rifie des propri{\'e}t{\'e}s
analogues aux 
pr{\'e}c{\'e}dentes~: pour  $x=(x_1,x_2) \in \wt X$ et $\theta  = (\theta
_1,\theta _2) \in \pa \hf 1 \times \pa \hf 2$
$$d\mu _{x} = e^{-(n_1-1)B_1(x_1,\theta _1)-(n_2-1)B_2(x_2,\theta _2)} d\theta_1\otimes d\theta _2~ .$$
\begin{remarque}
Remarquons que la mesure ci-dessus est diff\'erente de celle utlis\'ee dans les r\'ef\'erences \cite{Alb}, \cite{Quint} et \cite{Link}.
\end{remarque}
En effet,

{\it i)\/}~~Pour $O$ et $x\in \wt X$, $d\mu _O$ et $d\mu _x$ sont
absolument
continues, mais la densit{\'e} n'a plus la forme pr{\'e}c{\'e}dente,
elle vaut~: 
$${d\mu _x\over d\mu _O} = \exp\left(-\left[(n_1-1) B_1(x_1,\theta
_1)+(n_2-1) B_2(x_2,\theta _2)\right]\right)~.$$
On remarque que $\mu ^i_{x_i} = e^{-(n_i-1)B_i(x_i\theta _i)}d\theta
_i$ est la mesure de Patterson-Sullivan de $(\hf i,g^i_0)$.

\medskip
{\it ii)\/}~~Si $\gamma  = (\gamma _1,\gamma _2) \in \Isom (\hf 1,g^1_0)
\times \Isom(\hf 2,g^2_0)$ alors
$$\mu _{\gamma (x)} = \gamma _*(\mu _x)$$
car
\begin{align*}
\gamma _*(\mu _x)&= (\gamma _1,\gamma _2)_* (\mu ^1_{x_1} \otimes \mu
^2_{x_2}) = (\gamma _1)_*(\mu ^1_{x_1}) \otimes  (\gamma _2)_*(\mu
^2_{x_2})\\
&= \mu ^1_{\gamma _1(x_1)} \otimes \mu ^2_{\gamma _2(x_2)} = \mu
_{\gamma (x)}~.\end{align*}
De m{\^e}me, si $n_1 = n_2$, on v{\'e}rifie ais{\'e}ment que l'isom{\'e}trie
suppl{\'e}mentaire
$$\zeta (x_1,x_2) = (x_2,x_1)$$
satisfait cette contrainte.

Dans la suite nous travaillerons donc avec cette famille $\mu _x$ qui
est le produit des mesures de Patterson-Sullivan de chaque
facteur. Terminons en remarquant que si $B^{\alpha_i} _i$ d{\'e}signe la
fonction de Busemann de $(\hf i,\alpha_i ^2g^i_0)$, alors
$$\Ent(\alpha_i ^2g^i_0) B^{\alpha_i} _i(\cdot ,\cdot) = {1\over \alpha_i }
(n_i-1) \alpha_i  B^i(\cdot,\cdot)~;$$
de sorte que la famille $\mu
_x$ ne d{\'e}pend ni de $\alpha_1$, ni de $\alpha_2 $.

\subsection{Barycentre}

Nous construisons ici une application inverse de $x\mapsto \mu _x$,
\cad une application qui associe {\`a} la plupart des mesures sur $\pa_F
\wt X$ un point de $\wt X$ qui est son centre de masse ou
barycentre. La construction est analogue {\`a} celle de \cite{BCG1} et \cite{BCG2} {\`a}
l'utilisation pr{\`e}s de $\pa_F \wt X$ au lieu de $\pa \wt X$.

Soit $\vc$  une mesure bor{\'e}lienne positive non nulle  sur
$\pa_F \wt X$, on consid{\`e}re la fonction
$$x\in \wt X,~~ \bc_{\alpha_1 ,\alpha_2 }(x) = \int_{\pa_F \wt X} B_{\alpha_1
,\alpha_2 } (x,\theta ) \ d\vc (\theta )~.$$
On d{\'e}finit les mesures marginales sur $\pa \hf 1$ et $\pa \hf 2$ par~:

{\it i)\/}~~$\vc _1(A_1) = \vc (A_1 \times \pa \hf 2) = \pi _{1*}
(\vc)$, o{\`u} $A_1$ est un bor{\'e}lien de $\pa \hf 1$ et $\pi _1$  la
projection canonique de $\pa_F \wt X$ sur $\pa \hf 1$~; et de m{\^e}me,

{\it ii)\/}~~  $\vc _2(A_2) = \vc ( \pa \hf 1 \times A_2) = \pi _{2*}
(\vc)$, o{\`u} $A_2$ est un bor{\'e}lien de $\pa \hf 2$ et $\pi _2$ la
projection  de $\pa_F \wt X$ sur $\pa \hf 2$~.

\begin{proposition}
Si $\vc_1$ et $\vc_2$ sont des mesures non nulles
et sans atomes, pour tous $\alpha_1 $, $\alpha_2 $ strictement positifs, la
fonction $\bc_{\alpha_1 ,\alpha_2 }$ est $\ci$, strictement convexe sur
$\wt X$ et tend vers l'infini lorsque $x$ tend vers l'infini.
\end{proposition}

\begin{preuve}
Par d{\'e}finition de $\vc_1$, $\vc_2$ et $B_{\alpha_1 ,\alpha_2 }$, on
a~:
$$\bc_{\alpha_1 ,\alpha_2 }(x) = {1\over \sqrt{\alpha_1 ^2+\alpha_2 ^2}} \left(
\alpha_1 \int_{\pa \hf 1} B_1(x_1,\theta _1)\ d\vc_1(\theta _1) +
\alpha_2 \int_{\pa \hf 2} B_2(x_2,\theta _2)\ d\vc_2(\theta _2)
\right).$$
En effet,
\begin{align*}
\int_{\pa _F\wt X} B_1\big(x_1,\theta _1) \ d\vc(\theta _1,\theta _2)\big) &=
\int_{\pa_F \wt X} B_1\big(x_1,\pi _1(\theta _1,\theta _2)\big) \ d\vc
(\theta _1,\theta _2)\\
&= \int_{\pa \hf 1} B_1(x_1,\theta ) d(\pi _1*\vc)(\theta _1)\end{align*}
et de m{\^e}me avec l'autre terme. Alors, on applique les r{\'e}sultats de
\cite{Dou-Ear}, \cite{BCG1} et \cite{BCG2} qui montrent que $x_i \mapsto \int_{\pa \hf 1}
B_i(x_i,\theta _i) \ d\vc(\theta _i)$ est strictement convexe, pour
$i=1,2$, et tend vers l'infini lorsque $x_i$ tend vers l'infini dans
$\hf i$. On rappelle qu'une fonction est dite strictement convexe si
elle l'est en restriction {\`a} toute g{\'e}od{\'e}sique non constante. Il est
alors facile de v{\'e}rifier que $\bc_{\alpha_1 ,\alpha_2 }$ est strictement
convexe en restriction {\`a} toute g{\'e}od{\'e}sique non constante de $\wt X =
\hf1 \times \hf2$. Les autres conclusions de la proposition sont
{\'e}galement {\'e}videntes.
\end{preuve}

\begin{remarque}
L'hypoth{\`e}se sur la mesure $\vc$ est v{\'e}rifi{\'e}e, par
exemple, d{\`e}s que celle-ci est absolument continue par rapport {\`a} la
mesure de Lebesgue sur $\pa_F \wt X$. Par ailleurs, elle peut {\^e}tre
affaiblie (voir \cite{BCG1}).
\end{remarque}

\begin{corollaire}\label{barycentre}
Sous les m{\^e}mes hypoth{\`e}ses, la fonction
$\bc_{\alpha_1 ,\alpha_2 }$ admet un unique minimum sur $\wt X$ que nous
appellerons le barycentre de $\vc$, not{\'e} $\barr(\vc)$, qui ne d{\'e}pend
pas de $\alpha_1 $, $\alpha_2 $ ({\`a} condition qu'ils soient strictement
positifs). De plus $\barr(\vc) = (\barr_1(\vc_1), \barr_2(\vc_2))$, o{\`u}
$\barr_i(\vc_i)$ d{\'e}signe le barycentre de la mesure $\vc_i$ dans $\hf i$.
\end{corollaire}

\begin{preuve}
L'unicit{\'e} r{\'e}sulte de la stricte convexit{\'e} de
$\bc_{\alpha_1 ,\alpha_2 }$ 
$(\alpha_1 >0,~\alpha_2 >0)$ et du fait $\bc_{\alpha_1 ,\alpha_2 }(x) \build
\hfl_{x\to +\infty}\fin  +\infty$. Le point $x^* = (x^*_1,x^*_2)$ est
d{\'e}fini par l'{\'e}quation vectorielle,
$$\overrightarrow{\nabla \bc_{\alpha_1 ,\alpha_2 }}(x^*) = \overrightarrow
O$$
c'est-\`a-dire $\alpha_1 \int_{\pa \hf 1} \overrightarrow{\nabla_1
\bc_1}(x^*_1,\theta _1)\ d\vc_1(\theta _1) + \alpha_2  \int_{\pa \hf
2} \overrightarrow{\nabla_2 
\bc_2}(x^*_2,\theta _2)\ d\vc_2(\theta _2) =0$ (ici $\nabla_i$ d{\'e}signe
le gradient d'une fonction d{\'e}finie sur $\hf i$).

Si $\bar x_i = \barr_i(\vc_i)$ (voir \cite{BCG2}), alors
$$\int_{\pa \hf i} \nabla_i B_i (\bar x_i,\theta _i) \ d\vc(\theta _i)
= \overrightarrow O \hbox{~~pour~~} i=1,2~.$$
Par unicit{\'e} on a donc $\bar x_i = x^*_i$ ($i=1,2$), \cad
$$\barr(\vc_{\alpha_1 ,\alpha_2 }) = (\barr_1(\vc_1),
\barr_2(\vc_2))~.$$
\end{preuve}

\section{Le th{\'e}or{\`e}me principal}

Dans ce chapitre, nous nous proposons de prouver le th{\'e}or{\`e}me principal
sur l'entropie, analogue, dans cette situation, des r{\'e}sultats prouv{\'e}s
dans \cite{BCG1} et \cite{BCG2}. Nous donnons l'{\'e}nonc{\'e} et la preuve dans un cas
particulier afin d'{\'e}viter des lourdeurs dans les notations~; le cas
g{\'e}n{\'e}ral est rigoureusement identique.

Nous consid{\'e}rons, comme pr{\'e}c{\'e}demment, $\wt X = \hhf^{n_1} \times \hhf^{n_2}$,
o{\`u} $n_i \ge 3$  et nous munissons $X$ de la m{\'e}trique $g_0 =
g_{a_1,a_2}$, o{\`u} les nombres $a _i$ sont ceux calcul{\'e}s dans la
proposition 1.4. La m{\'e}trique $g_0$ minimise l'entropie normalis{\'e}e, sur
$X$, parmi les m{\'e}triques $g_{\alpha _1,\alpha _2}$  (voir la
proposition \ref{min}). Soit $\Gamma $ un sous-groupe du groupe d'isom{\'e}tries de
$(\wt X,g_0)$ tel que $X = \wt X / \Gamma $ est une vari{\'e}t{\'e} compacte
($\Gamma $ est un r{\'e}seau co-compact et sans torsion).

\begin{theoreme}\label{principal}
Soit $(Y,g)$ une vari{\'e}t{\'e}
riemannienne compacte de 
dimension $n=n_1+n_2$ et $f : Y\to X$ une application continue, alors

\begin{itemize}
\item[i)] $(\Ent (Y,g))^n \vol (Y,g) \ge |\deg f| \Ent (X,g_0)^n
\vol (X,g_0)$~;

\item[ii)] l'{\'e}galit{\'e}, dans l'in{\'e}galit{\'e} ci-dessus, a
lieu si, et 
seulement si, $f$
est homotope {\`a} un rev{\^e}tement riemannien.
\end{itemize}
\end{theoreme}

\begin{remarque}
Ce r{\'e}sultat est vrai dans la version g{\'e}n{\'e}rale
donn{\'e}e en introduction. Sa preuve est analogue {\`a} celle de \cite{BCG1}.
\end{remarque}

\vskip5pt
{\it  Preuve de l'in{\'e}galit{\'e} i).}

Nous donnons une preuve
inspir{\'e}e de la 
technique d{\'e}velopp{\'e}e dans \cite{BCG3}. On note $\Gamma  = \pi _1(Y)$, le
groupe fondamental de $Y$, $\wt Y$ le rev{\^e}tement universel de
$Y$. L'application continue $f : Y\to X$ induit un morphisme $\rho  :
\Gamma  \to \Gamma _0  = \pi _1(X)$. On appelle $\mu _O$ la mesure
(canonique) $d\theta _1 \otimes d\theta _2$ sur $\pa_F \wt X$ et
$\mc(\pa_F \wt X)$ l'espace des mesures de Radon positives sur
$\pa_F \wt X$. Soit $O\in \wt Y$ un point fix{\'e} (une origine)~;
consid{\'e}rons l'application
\begin{align*}
\wt Y &\la \mc(\pa_F \wt X)\\
y &\longmapsto \mu _{y,\varepsilon } = \summ_{\gamma \in \Gamma }
e^{-\Ent(Y,g) (1+\varepsilon ) d_{\wt Y} (y,\gamma (O))} \rho (\gamma
)_*(\mu _O) \end{align*}
o{\`u} $\varepsilon >0$. Cette application est {\'e}quivariante; en effet,
pour tout $\alpha \in \Gamma $,
$$\mu _{\alpha  (y),\varepsilon } = \rho (\alpha  )_* (\mu _{y,\varepsilon })~.$$
On d{\'e}finit alors
\begin{align*}
\wt F_\varepsilon  : \wt Y &\la \wt X\\
y&\longmapsto  \barr(\mu _{y,\varepsilon })~.\\\end{align*}
Notons que chaque mesure $\mu_{y,\varepsilon }$ est sans atome.
La notion de barycentre {\'e}tant ind{\'e}pendante des coefficients $\alpha _i$
servant {\`a} d{\'e}finir la m{\'e}trique, nous utiliserons, pour simplifier, la
m{\'e}trique $\bar g_0 = g^1_0 \oplus g^2_0$ (voir le corollaire \ref{barycentre}).

Alors, par {\'e}quivariance, $\wt F_\varepsilon $ donne une famille d'applications
$$F_\varepsilon  : Y\to X~.$$
Par ailleurs, le barycentre sur $\wt X$ se d{\'e}compose (\cf
corollaire \ref{barycentre}) et donc {\'e}galement la fonction $\wt F_\varepsilon $
\begin{align*}
\wt F_\varepsilon  : \wt Y &\la \wt X = \hhf^{n_1} \times \hhf^{n_2}\\
y &\longmapsto (\wt F_{1,\varepsilon } (y),\wt F_{2,\varepsilon }(y))\\\end{align*}
o{\`u} $\wt F_{i,\varepsilon }(y) = \barr_i (\pi _i * (\mu _{y,\varepsilon
}))$.

Nous notons $F_\varepsilon  $, $F_{1,\varepsilon }$ et $F_{2,\varepsilon }$
les applications correspondantes de $Y$ dans $X$.

\begin{lemme}
Pour tout $\varepsilon >0$, les fonctions
$F_\varepsilon $, $F_{1,\varepsilon }$ et $F_{2,\varepsilon }$ sont
lipschitziennes.
\end{lemme}

\begin{preuve} 

{\it i)\/} La s{\'e}rie $\summ_{\gamma \in \Gamma } e^{-\Ent (Y,g)
(1+\varepsilon ) d_{\wt Y} (y,\gamma (O))}$ converge d{\`e}s que
$\varepsilon >0$~; en effet, puisque $Y$ est compacte, elle est
comparable {\`a} l'int{\'e}grale 
$$I_\varepsilon  = \int_{\wt Y} e^{-\Ent(Y,g) (1+\varepsilon ) d_{\wt
Y} (y,z)} dv_g(z)$$
qui converge. Pour v{\'e}rifier ce dernier point il suffit d'{\'e}crire
$I_\varepsilon $ en coordonn{\'e}es polaires et d'appliquer la d{\'e}finition
de $\Ent(Y,g)$.

{\it ii)\/} La fonction $\wt F_\varepsilon $ est d{\'e}finie par l'{\'e}quation
implicite
$$\int_{\pa_F \wt X} \overrightarrow{\nabla {\ol B}_0}_{|(F_\varepsilon
(y),\theta )} (\summ_{\gamma \in \Gamma } e^{-\Ent (Y,g)
(1+\varepsilon ) d_{\wt Y} (y,\gamma (O))} \rho (\gamma )_* (d\mu _O)) =
\overrightarrow 0$$
o{\`u} $\ol B_0$ d{\'e}signe la fonction de Busemann de la m{\'e}trique ${\ol g}_0 =
g_{1,1}$ sur $\wt X$. Ici  $\theta  = (\theta _1,\theta
_2)$. L'{\'e}quation ci-dessus peut se r{\'e}crire en
$$0 = L(x,y) = \int_{\pa_F \wt X} \summ_{\gamma  \in \Gamma } 
\overrightarrow{\nabla {\ol B}_0}_{(x,\rho (\gamma )\theta )} e^{-\Ent (Y,g)
(1+\varepsilon ) d_{\wt Y} (y,\gamma (O))} d\mu _O $$
o{\`u} $x = (x_1,x_2) \in
\wt X$. La fonction $L$ est $\ci$ en $x$, lipschitzienne en $y$ et
chaque diff{\'e}rentielle partielle en $x$ est lipschitzienne en $y$. Dans
cette situation, on peut utiliser le th{\'e}or{\`e}me des fonctions implicites
pour conclure que $\wt F_\varepsilon $ est lipschitzienne en $y$
(voir~\cite{Dieu}). Remarquons que la condition sur la diff{\'e}rentielle
partielle en $x$ qui est n{\'e}cessaire pour appliquer le th{\'e}or{\`e}me des
fonctions implicites est exactement celle qui prouve l'existence du
barycentre (voir le chapitre pr{\'e}c{\'e}dent), \cad la stricte
convexit{\'e} de 
$\bc_{1,1}$.
\end{preuve}
\begin{lemme}\label{jacobien}
Pour $\varepsilon >0$ et pour tout $y \in \wt Y$, on a
$$|\Jac {\wt F_\varepsilon} (y)| \le (1+\varepsilon )^n \Big({\Ent
(Y,g)\over \Ent(X,g_0)} \Big)^n$$
o{\`u} $n = \dim Y = \dim X = n_1 +n_2$ et le Jacobien est calcul{\'e} {\`a}
l'aide des m{\'e}triques $g$ sur $Y$ et $g_0$ sur $X$.
\end{lemme}

\begin{preuve}
Comme nous l'avons remarqu{\'e} dans le chapitre pr{\'e}c{\'e}dent la notion de
barycentre, et donc la d{\'e}finition de l'application $\wt F$, ne d{\'e}pend pas
de $\alpha _1$, $\alpha _2$. Nous pouvons donc
utiliser sur $\wt X$ la m{\'e}trique $\bar g_0 = g^1_0 \oplus g^2_0$ (on
rappelle que $g^i_0$ d{\'e}signe ici la m{\'e}trique de courbure constante
{\'e}gale {\`a} $-1$ sur
$\hhf^{n_i}$). Rappelons {\'e}galement la notation $g_0 = a_1 g^1_0
\oplus a_2 g^2_0$ o{\`u} $a_i$ sont les valeurs calcul{\'e}es dans la section
pr{\'e}c{\'e}dente, telles que $g_0$ minimise l'entropie
normalis{\'e}e parmi les 
m{\'e}triques $g_{\alpha _1,\alpha _2}$. Nous noterons
$\overline{\det}(D{\wt F_\varepsilon}  (y))$ le d{\'e}terminant de la
diff{\'e}rentielle de ${\wt F_\varepsilon} $ en $y$ calcul{\'e} {\`a} l'aide des
m{\'e}triques $g$ sur $\wt Y$ et $\overline{g_0}$ sur $\wt X$~; par
ailleurs  $\Jac
{\wt F_\varepsilon} (y) = a^{n_1}_1 a^{n_2}_2 \overline{\det}
(D{\wt F_\varepsilon}  (y))$ est le d{\'e}terminant de $D{\wt F_\varepsilon} (y)$ calcul{\'e}
{\`a} l'aide des m{\'e}triques $g$ sur $\wt Y$ et $g_0$ sur $\wt X$. Notons
que $g_0$ est normalis{\'e}e par $a^{n_1}_1 a^{n_2}_2 = 1$, de sorte que
$\Jac {\wt F_\varepsilon} (y) = a^{n_1}_1 a^{n_2}_2 \ol{\det}
(D{\wt F_\varepsilon} (y)) = \ol{\det} (D{\wt F_\varepsilon} (y))$. Nous
distinguerons toutefois les deux expressions afin d'{\'e}viter les
confusions entre les m{\'e}triques $g_0$ et $\bar g_0$.

\medskip
 {\it Estimation de $\overline{\det}(D{\wt F_\varepsilon} (y))$}; ici, tous les
calculs se font {\`a} l'aide de la m{\'e}trique $\bar g_0$ sur $\wt X$.
Rappelons que nous d{\'e}signons par $B_i$ les fonctions de Busemann sur
$\hhf^{n_i}$ muni de la m{\'e}trique $g^i_0$. Comme dans \cite{BCG3}, page 155,
nous posons
\begin{align*}
k_{y,\varepsilon }(v,v) &= \int_{\pa_F \wt X} Dd\ol B_{0}
{}_{|({\wt F_\varepsilon} (y),\theta )}(v,v)\  d\mu _{y,\varepsilon }
(\theta ) = \bar g_0 (K_{y,\varepsilon } (v),v)\\
h_{y,\varepsilon }(v,v) &= \int_{\pa_F \wt X} \Big( d\ol B_{0}{}_{|({\wt F_\varepsilon} (y),\theta
)} (v)\Big)^2 \ d\mu _{y,\varepsilon }(\theta )\\
&= \bar g_0 (H_{y,\varepsilon }(v),v)\\\end{align*}
o{\`u} $v \in T_{{\wt F_\varepsilon} (y)} \wt X$. Ici, comme dans la section
pr{\'e}c{\'e}dente,
$$\ol B_{0}(x,\theta ) = {1\over \sqrt 2} \big(B_1(x_1,\theta _1) +
B_2(x_2,\theta _2)\big)~.$$
Enfin,
\begin{align*}
h'_{y,\varepsilon }(u,u) 
&= \summ_{\gamma \in \Gamma } \int_{\pa_F\wt X} \langle \nabla d_{\wt
Y| (y,\gamma (O))},u\rangle^2 e^{-\Ent (Y,g) (1+\varepsilon )
d_{\wt Y} (y,\gamma (O))}\ d(\rho (\gamma )_*\mu_O )(\theta )\\
&= \big(\mu_O (\pa_F\wt X)\big) \summ_{\gamma \in \Gamma }\langle
\nabla d_{\wt Y| (y,\gamma (O))},u\rangle^2 e^{-\Ent (Y,g) (1+\varepsilon )
d_{\wt Y} (y,\gamma (O))}\\
&= g(H'_{y,\varepsilon } u,u)\end{align*}
pour $u\in T_y \wt Y$. Nous utiliserons les m{\^e}mes notations pour les
formes bilin{\'e}aires associ{\'e}es. En diff{\'e}renciant
l'{\'e}quation implicite 
qui d{\'e}finit ${\wt F_\varepsilon} $, nous obtenons, pour $u \in T_y \wt Y$ et
$v \in T_{{\wt F_\varepsilon} (y)} \wt X$,
\begin{multline*}
k_{y,\varepsilon } \big(v,D {\wt F_\varepsilon} (y)(u)\big) =
(1+\varepsilon )\Ent (Y,g) \summ_{\gamma \in \Gamma } \int_{\pa_F \wt
X} g\Big(\nabla d_{\wt Y|(y,\gamma (O))},u\Big) \bar g_0
\Big(\nabla \ol B_{0|({\wt F_\varepsilon} (y),\theta )},v\Big)\hfill\\
\hfill \cdotp e^{-\Ent
(Y,g) (1+\varepsilon ) d_{\wt Y} (y,\gamma (O))} d(\rho (\gamma )_*\mu_O
)(\theta )
\end{multline*}
et, en utilisant l'in{\'e}galit{\'e} de Cauchy-Schwarz,
$$k_{y,\varepsilon }\big(D {\wt F_\varepsilon}  (y) (u),v\big) \le
(1+\varepsilon ) \Ent (Y,g) \big( h_{y,\varepsilon } (v,v)\big)^{1/2}
\big( h'_{y,\varepsilon } (u,u)\big)^{1/2}~.\leqno(2.5)$$
Un lemme {\'e}l{\'e}mentaire d'alg{\`e}bre lin{\'e}aire (\cf \cite{BCG2},
lemme 5.4) donne, 
{\`a} partir de (2.5),
$$\det(K_{y,\varepsilon }) \overline{\det}(D{\wt F_\varepsilon} (y)) \le
\big((1+\varepsilon ) \Ent (Y,g)\big)^{n_1+n_2} \big(\det H_{y,\varepsilon}
\big)^{1/2} \big(\det H'_{y,\varepsilon}
\big)^{1/2} ~.$$

Nous devons maintenant  remarquer que la notion de barycentre ne
change pas lorsque l'on multiplie une mesure par un  nombre
strictement positif, de sorte que
$${\wt F_\varepsilon} (y) = \barr (\mu _{y,\varepsilon }) = \barr \Big({\mu
_{y,\varepsilon }\over \mu _{y,\varepsilon } (\pa_F \wt X)}\Big)~.$$
On peut donc supposer que la famille de mesures que l'on consid{\`e}re est
normalis{\'e}e (de masse totale {\'e}gale {\`a} 1 pour tout $y \in \wt Y $ et
$\varepsilon  >0$). La trace d'une forme quadratique $\varphi $
(calcul{\'e}e dans une base orthonorm{\'e}e par rapport {\`a} une structure
euclidienne $g$) {\'e}tant not{\'e}e $\tr_g \varphi $, en injectant dans la
d{\'e}finition de $h'_{y,\varepsilon }$ le fait que $\|\nabla d_{\wt
Y}\|_g =1$, nous obtenons
$$\tr (H'_{y,\varepsilon }) = \tr_g (h'_{y,\varepsilon }) = 1~,$$
d'o{\`u}
$$\big(\det H'_{y,\varepsilon }\big)^{1/2} \le \Big({1 \over
\sqrt{n_1+n_2}}\Big)^{n_1+n_2}~.$$
Maintenant la d{\'e}finition de $h_{y,\varepsilon }$ (et $H_{y,\varepsilon
}$) montre que
$$H_{y,\varepsilon } = 
\begin{pmatrix}
H_1  &*\\
*  &H_2
\end{pmatrix}
$$
o{\`u} $H_i$ d{\'e}signe la restriction de $H_{y,\varepsilon }$ {\`a}
$\hhf^{n_i}$~; plus pr{\'e}cis{\'e}ment, pour $i=1,2$ et $v _i \in
T_{F_{i,\varepsilon }(y)} \hhf^{n_i}$
\begin{align*}
g^i_0 (H_iv_i,v_i) &= \int_{\pa_F \wt X} \Big(
d\ol B_{0|({\wt F_\varepsilon} (y),\theta )} (v_i)\Big)^2\ d\mu
_{y,\varepsilon } (\theta )\\
&= \int_{\pa_F\wt X} {1\over 2}\Big(
dB_{i|({\wt F_{i,\varepsilon}} (y),\theta_i )} (v_i)\Big)^2\ d\mu
_{y,\varepsilon } (\theta )\\
&= {1\over 2} \int_{\pa \hhf^{n_i}}\Big(
dB_{i|({\wt F_{i,\varepsilon}} (y),\theta_i )} (v_i)\Big)^2\ d\big( \pi
_{i*} (\mu _{y,\varepsilon })\big) (\theta_i )~.\\  \end{align*}
Remarquons que, puisque $\|dB_i\| _{g^i_0} = 1$ et $(\pi_i)_*\mu_{y,\varepsilon}$ est une probabilit\'e, nous avons
$$\tr (2H_i) = 1~.$$
De m{\^e}me,
$$K_{y,\varepsilon } =
\begin{pmatrix} 
K_1  &0\\
0  &K_2
\end{pmatrix}
$$
avec, pour $i=1,2$ et $v_i = T_{{\wt F_{i,\varepsilon }}(y)} \hhf^{n_i}$,
$$g^i_0 (K_iv_i,v_i) = {1\over \sqrt 2} \int_{\pa\hhf^{n_i}}
DdB_{i| ({\wt F_{i,\varepsilon }}(y),\theta _i)} (v_i,v_i)\ d\big(\pi
_{i*} (\mu _{y,\varepsilon })\big)(\theta _i)~.$$

\begin{lemme} \label{2.6}
Avec les notations pr{\'e}c{\'e}dentes, nous avons

\begin{itemize}
\item[i)] $\det(K_{y,\varepsilon }) = \det (K_1) \det (K_2)~;$

\item[ii)]  $\det(H_{y,\varepsilon }) \le \det (H_1) \det (H_2)$.
\end{itemize}
\end{lemme}

L'{\'e}galit{\'e} {\it i)\/} est {\'e}vidente et l'in{\'e}galit{\'e}
{\it ii)\/} est 
classique (voir \cite{Bec-Bel}, p.~63) pour les matrices sym{\'e}triques.

Par ailleurs, sur les espaces hyperboliques $\hhf^{n_i}$, la relation
suivante est v{\'e}rifi{\'e}e (voir \cite{BCG1}, p.~751), pour $i=1,2$, 
$$DdB_i = g^i_0 - dB_i \otimes dB_i$$
qui se traduit en
$$K_i = {1\over \sqrt 2} (I_i - 2H_i)$$
o{\`u} $I_i$ d{\'e}signe l'identit{\'e} de $T_{F_{i,\varepsilon }(y)}
\hhf^{n_i}$. En regroupant ces remarques, nous obtenons, {\`a} partir de
2.5,
$$\overline{\det}(D{\wt {\wt F_\varepsilon}} (y))
\le \bigg({(1+\varepsilon )\Ent(Y,g )\over \sqrt{n_1+n_2}}
\bigg)^{n_1+n_2}
{(\det 2H_1)^{1/2} \over \det (I-2H_1)}~ {(\det 2H_2)^{1/2} \over \det
(I-2H_2)} ~.$$

Alors, un lemme alg{\'e}brique donne (\cf \cite{BCG1}, appendice B),
$${(\det 2H_i)^{1/2} \over \det (I-2H_i)} \le \Big({\sqrt{n_i}\over
n_i-1}\Big)^{n_i}$$
l'{\'e}galit{\'e} n'ayant lieu que si, et seulement si, $2H_i = {1\over n_i}
I_i$ (on rappelle que $\tr (2H_i) = 1$).

\medskip


En regroupant ces in{\'e}galit{\'e}s, il vient
\begin{align*}
\Jac {\wt {\wt F_\varepsilon}} (y)
&= a^{n_1}_1 a^{n_2}_2 \big( \overline{\det}(D{\wt {\wt F_\varepsilon}}
(y)\big)\\
&\le \big((1+\varepsilon ) \Ent (Y,g)\big)^{n_1+n_2} {a^{n_1}_1
a^{n_2}_2\over (\sqrt{n_1+n_2})^{n_1+n_2}}
\Big({\sqrt{n_1}\over n_1-1}\Big)^{n_1}~\Big({\sqrt{n_2}\over
n_2-1}\Big)^{n_2}\\
&= \bigg((1+\varepsilon ) {\Ent (Y,g)\over \Ent(X,g_0)} \bigg)^n\\\end{align*}
d'apr{\`e}s la proposition \ref{min}. Ce qui prouve le lemme \ref{jacobien}.

L'in{\'e}galit{\'e} {\it i)\/} du th{\'e}or{\`e}me \ref{principal} s'en d{\'e}duit
par int{\'e}gration et 
passage {\`a} la limite en $\varepsilon =0$.

Le cas d'{\'e}galit{\'e} sera trait{\'e}, dans un cadre plus
g{\'e}n{\'e}ral, dans le 
paragraphe suivant.

\vskip5pt
{\it  Remarques sur le cas g{\'e}n{\'e}ral\pointir} Si $\wt X = \wt X_1 \times \cdots
\times \wt X_p$ et $\bar g_0 = g_1 \oplus \cdots \oplus g_p$ o{\`u} $(\wt
X_k,g_k)$ est un espace sym{\'e}trique de courbure strictement n{\'e}gatif et
de dimension $n_k$, on munit $X$ de la m{\'e}trique $g_0 = a_1g_1 \oplus
\cdots \oplus a_p g_p$, o{\`u} les nombres $a_i$ sont ceux calcul{\'e}s dans
la proposition \ref{mingeneral}. La m{\'e}trique $g_0$ minimise l'entropie normalis{\'e}e
parmi les m{\'e}triques $g_\alpha $ (voir la proposition \ref{mingeneral}).

Alors, comme ci-dessus, on pose $\bar g_0 = g_1 \oplus \cdots \oplus
g_p$. On suppose de plus que la courbure sectionnelle de $(\wt
X_k,g_k)$ est normalis{\'e}e de sorte qu'elle soit {\'e}gale {\`a} $-1$ si $(\wt
X_k,g_k)$ est hyperbolique r{\'e}elle et comprise entre $-4$ et $-1$ dans
les autres cas. Le calcul de l'entropie d'une telle m{\'e}trique est donn{\'e}
dans \cite{BCG1}, p.~740.

Pour $x = (x_1\ld x_p) \in \wt X$ et $\theta  = (\theta _1\ld \theta
_p) \in \partial_F \wt X$
 ($\partial_F \wt X = \partial \wt X_1 \times \cdots \times \partial
\wt X_p$), la fonction de Busemann de $(\wt X,\bar g_0)$ est
$$\ol B_0(x,\theta ) = {1 \over \sqrt{p}} \Big(B_1(x_1,\theta _1) + \cdots +
B_p(x_p,\theta _p)\Big)$$
et on a les d{\'e}compositions
$$
H_{y,\varepsilon } = \begin{pmatrix}
H_1 &*&*&*\\
* &H_2 &*&*\\
* &*&\ddots &*\\
*&*&*&H_p\\
\end{pmatrix}
,\qquad 
K_{y,\varepsilon} = \begin{pmatrix}
 K_1\kern1pt &0&0&0\\
0 &K_2\kern1pt &0&0\\
0 &0&\ddots &0\\
0&0&0&K_p\kern1pt \\
\end{pmatrix}
$$
avec $\tr (pH_k) =1$ pour $k = 1,2\ld p$.

La relation qui lie $K_i$ et $H_i$ d{\'e}pend du type d'espace consid{\'e}r{\'e}
(hyperbolique r{\'e}el, complexe, quaternionien ou de Cayley) et est
d{\'e}crite dans \cite{BCG1}, p.~751. On peut v{\'e}rifier ais{\'e}ment que
$$\tr         (\sqrt p K_k) = E_k = \hbox{~entropie de~} (\wt X_k,g_k)$$
pour $k = 1\ld p$. 

Dans l'appendice $B$ de \cite{BCG1} nous montrons que
$${\det (pH_k)^{1/2}\over \det(\sqrt p K_k)} \le
\biggl({\sqrt{n_k}\over E_k}\biggr)^{n_k}~.$$
On conclut, alors, gr{\^a}ce {\`a} la proposition \ref{mingeneral}, comme ci-dessus.
\end{preuve}
\section{Le volume des repr{\'e}sentations}

Nous donnons dans ce paragraphe une application de la technique
introduite dans \cite{BCG3} aux repr{\'e}sentations du groupe fondamental d'une
vari{\'e}t{\'e} compacte.

Dans ce qui suit $\wt X$ est un produit fini d'espaces sym{\'e}triques
simplement connexe de courbure strictement n{\'e}gative. Chaque facteur
est suppos{\'e} de dimension sup{\'e}rieure ou {\'e}gale {\`a} 3. On munit $\wt X$  de
la m{\'e}trique $g_0$ d{\'e}crite dans la proposition \ref{min}, \cad celle qui
r{\'e}alise l'entropie minimale pour tous les quotients compacts de $\wt
X$. Par ailleurs, $(Y,g)$ est une vari{\'e}t{\'e} riemannienne compacte dont
le groupe fondamental est not{\'e} $\Gamma $. On consid{\`e}re
$$\rho  : \Gamma  \la \Isom (\wt X,g_0)$$
une repr{\'e}sentation. Il existe toujours des applications {\'e}quivariantes
$f : \wt Y \to \wt X$ car $\wt X$ est contractile (dans la suite nous
donnerons un exemple explicite d'une telle application). Elle v{\'e}rifie
donc
$$\forall \varphi \in \Gamma , ~ \forall y \in \wt Y,~~ f(\gamma (y)) =
\rho (\gamma )\ f(y)~.$$
On peut toujours la supposer $C^1$, quitte {\`a} la r{\'e}gulariser. Si on note
$\omega _0$ la forme volume de $(\wt X,g_0)$ alors,

\begin{definition}
On appelle volume de la repr{\'e}sentation $\rho $, le
nombre
$$\vol(\rho ) = \int_Y f^* (\omega _0)~.$$
\end{definition}

\begin{remarques}

\begin{itemize}
\item[i)] La d{\'e}finition ci-dessus a un sens car, $f$ {\'e}tant $C^1$,
$f^*(\omega _0)$ est une forme continue sur $\wt Y$ qui de plus est
invariante par $\Gamma $. Par ailleurs, il est imm{\'e}diat de v{\'e}rifier
que $\vol(\rho )$ ne d{\'e}pend pas du choix de l'application {\'e}quivariante
$f$.

\item[ii)] Il faut interpr{\'e}ter $\vol(\rho )$ comme l'analogue
de la quantit{\'e} $|\deg f|\vol (X)$ du th{\'e}or{\`e}me \ref{principal}. En effet, lorsque
$\rho (\Gamma )$ est discret et cocompact, agissant sans point fixe,
nous nous trouvons dans la situation du paragraphe 3 o{\`u} $X=\wt X /\rho
(\Gamma )$ et $\vol(X) |\deg f| = \vol (\rho )$ par d{\'e}finition du
degr{\'e} de l'application $f$.
\end{itemize}
\end{remarques}

Nous prouvons donc un th{\'e}or{\`e}me analogue~:

\begin{theoreme}
\label{theo-princ}
Avec les notations ci-dessus~:
\begin{itemize}
\item[i)] $\vol(\rho ) \le \Bigl(\dis{\Ent(Y,g)\over \Ent(\wt
  X,g_0)}\Bigr)^n \vol (Y,g)$.

\item[ii)] L'{\'e}galit{\'e}, dans l'in{\'e}galit{\'e} ci-dessus a lieu si, et
seulement si, la repr{\'e}sentation $\rho $ est injective, $X  = \wt
X/\rho (\Gamma )$ est une vari{\'e}t{\'e} compacte et $(Y,g)$ est homoth{\'e}tique
{\`a} $(X,g_0)$.
\end{itemize}
\end{theoreme}

\begin{remarques}
\begin{itemize}
\item[i)] Ce r{\'e}sultat est un premier pas dans la compr{\'e}hension
des repr{\'e}sentations des groupes fondamentaux de vari{\'e}t{\'e}s compactes
dans des groupes d'isom{\'e}tries d'espaces sym{\'e}triques de type non
compact.

\item[ii)] Les exemples de telles repr{\'e}sentations sont rares
et nous discuterons ce point plus loin dans le texte. Plus rares
encore sont les exemples dont le volume est non nul.

\item[iii)] Seul le cas de dimension 2, o{\`u} notre m{\'e}thode ne
s'applique pas, est compl{\`e}tement compris (\cf \cite{Gol}). En particulier,
le th{\'e}or{\`e}me \ref{theo-princ} est une g{\'e}n{\'e}ralisation de la c{\'e}l{\`e}bre in{\'e}galit{\'e} de
Milnor-Wood (\cf \cite{Mil}, \cite{Woo} et \cite{Rez1}).
\end{itemize}
\end{remarques}

\begin{preuve} L'in{\'e}galit{\'e} est {\'e}l{\'e}mentaire et sa preuve est celle du th{\'e}or{\`e}me
\ref{principal}, {\it i)}. Le cas d'{\'e}galit{\'e} par contre est beaucoup plus difficile
car nous ne poss{\'e}dons pas de quotient compact de $\wt X$ ($\wt X/\rho
(\Gamma )$ n'est m{\^e}me pas un espace s{\'e}par{\'e}, en g{\'e}n{\'e}ral) sur lequel s'appuyer
afin d'utiliser la th{\'e}orie du degr{\'e} (voir la preuve du cas d'{\'e}galit{\'e}
de \cite{BCG1}).

Afin de traiter ce cas d'{\'e}galit{\'e} difficile nous devons consid{\'e}rer une
autre application {\'e}quivariante que celle introduite dans le paragraphe
3. Soit $f$ une premi{\`e}re application continue et $\rho $-{\'e}quivariante,
$$f : \wt Y \la \wt X~,$$
par exemple, nous pouvons prendre
comme pr{\'e}c{\'e}demment
$${\tilde f}(y) = \barr\Bigl( \summ_{\gamma \in \Gamma } e^{-\Ent(Y,g)
(1+\varepsilon )\ d(y,\gamma (O))} \rho (\gamma )_*d\mu \Bigr)$$
les notations {\'e}tant, ici, celles du paragraphe 2.

On rappelle que si $\theta  \in \pa_F\wt X$ et $z \in \wt X$,
$P_0(z,\theta )$ d{\'e}signe le noyau de Poisson de $\wt X$, normalis{\'e}  en
une origine $O_0\in \wt X$ de sorte que
$$P_0(O_0, \cdot) \equiv 1~.$$
Nous construisons une autre application, comme dans \cite{BCG1}, d{\'e}finie,
pour tout\break 
$c> \Ent(Y,g)$, par
$${\wt F}_c(y) = \barr \biggl(\Bigl(\int_{\wt Y} e^{-cd(y,z)} P_0\big(
{\tilde f}(z),\theta \big) \ dv_g(z)\Bigr) d\theta \biggl)~.$$

La preuve de l'in{\'e}galit{\'e} {\it i)\/} du th{\'e}or{\`e}me \ref{theo-princ} est rigoureusement
identique {\`a} celle donn{\'e}e dans le paragraphe 3. Nous ne la reproduirons
donc pas. Notons qu'elle peut \^etre faite \`a l'aide de la fonction $\tilde f$ d\'efinie e qu'il n'est pas n\'ecessaire d'utiliser la fonction ${\wt F_c}$; cette derni\`ere est toutefois beaucoup plus ais\'ee \`a manipuler dans la preuve du cas d'\'egalit\'e; elle est, par exemple plus r\'eguli\`ere que $f$.

Posons comme dans \cite{BCG1}, pour $c>\Ent(Y,g)$
$$\psi (c,y,\theta ) = \int_{\wt Y} e^{-cd(y,z)} P_0(f(z),\theta )\
dv_g(z)$$
et
$$\Phi(c,y,\theta ) = {\psi (c,y,\theta )\over \int_{\pa_F \wt X} \psi
(c,y,\theta )\ d\theta }$$
qui est de norme $L^1(\pa_F \wt X,d\theta )$ {\'e}gale {\`a} 1.

\begin{lemme}
L'application $(c,y) \mapsto \Phi (c,y,\cdot)$ est de
classe $C^1$ de l'intervalle $\ ]\!\Ent (Y,g),+\infty[\times \wt Y$ dans $L^1(\pa_F \wt
X)$.
\end{lemme}

\begin{preuve} |Il n'est pas possible de montrer le lemme ci-dessus par simple application du th\'eor\`eme de d\'erivation sous le signe somme. Toutefois, dans \cite{BCG1}, nous prouvons, comme corollaire du
th{\'e}or{\`e}me de convergence domin{\'e}e, que $y \mapsto \Phi (c,y,\cdot)$ est
de classe $C^1$ ({\`a} $c$ fix{\'e}) et, si $u \in T_y \wt Y$, sa
diff{\'e}rentielle est donn{\'e}e par
$$(u\cdot \psi )(c,y,\theta ) = -c \int_{\wt Y} e^{-cd(y,z)} (u\cdot
d)(y,z) P_0({\tilde f}(z),\theta )\ dv_g(z)$$
la continuit{\'e} en $c$ de cette quantit{\'e} est {\'e}vidente en remarquant que
$|u\cdot d| \le \|u\|_g$, que $P_0$ est strictement positif et que,
pour $y$ et $z$ fix{\'e}s, $c\mapsto e^{-cd(y,z)}$ est d{\'e}croissante en
$c$~; ceci permet d'appliquer une nouvelle fois le th{\'e}or{\`e}me de
convergence domin{\'e}e.

De  m{\^e}me, pour $y$ et $\theta $ fix{\'e}s, on peut appliquer le th{\'e}or{\`e}me de
d{\'e}rivation sous le signe somme afin de montrer la diff{\'e}rentiabilit{\'e} en
$c$ ({\`a} $y$ et $\theta $ fix{\'e}). En effet,
$$0\le d(y,z) e^{-cd(y,z)} \le e^{-c'd(y,z)}$$
pour tout $c'<c$. Ceci montre que
$${\pa \psi \over \pa c} (c,y,\theta ) = - \int_{\wt Y} d(y,z)
e^{-cd(y,z)} P_0\big({\tilde f}(z),\theta \big)\ dv_g(z)$$
existe et, encore gr{\^a}ce au th{\'e}or{\`e}me de convergence domin{\'e}e, est
continue en $(c,y)$. Ceci prouve le lemme ci-dessus. On remarque que le m{\^e}me
type d'argument que ceux utilis{\'e}s dans \cite{BCG1} montrent que ${\pa \psi
\over \pa c}$ est de classe $C^1$ comme fonction de $y$ {\`a} valeurs
dans $L^1(\pa_F \wt X)$.

De m{\^e}me $\psi $ est de classe $\ci$ en $c$ et chaque d{\'e}riv{\'e}e en $c$
est de classe $C^1$ en $y$ comme fonction de $\wt Y$ {\`a} valeurs dans
$L^1(\pa_F \wt X)$. L'assertion du lemme concernant $\Phi $ s'en d{\'e}duit.
\end{preuve}
\begin{lemme}
L'application 
\begin{align*}
{\wt F} : ]\!\Ent (Y,g),+\infty[ \times \wt Y &\la \wt X\\
(c,y) &\longmapsto {\wt F_c}(y)\end{align*}
est de classe $C^1$.
\end{lemme}

\begin{preuve} Il s'agit d'une simple application du th{\'e}or{\`e}me
des fonctions implicites (voir \cite{BCG1}). Rappelons la preuve de ce
fait. Soit $\{e_i(x)\}_{i=1\ld n}$ une base orthonorm{\'e}e de $T_x\wt X$
d{\'e}pendant de mani{\`e}re $\ci$ de $x \in \wt X$. D{\'e}finissons les fonctions
$$G_i(x,c,y) = \int_{\pa_F \wt X} d\ol B_0(x,\theta )(e_i(x)) \Phi
(c,y,\theta )\ d\theta $$
(on rappelle que $\ol B_0(x,\theta )$ d{\'e}signe la fonction de Busemann de
$(\wt X,g_0)$ normalis{\'e}e en $O_0$ et $d\theta $ la mesure canonique de
$\pa_F \wt X$), et
\begin{align*}
G : \wt X \times ]\!\Ent (Y,g),+\infty[\times \wt Y &\la \rrf^n\\
(x,c,y) &\longmapsto \big(G_1(x,c,y)\ld G_n(x,c,y)\big)~.\end{align*}
Alors, la fonction $\wt F$ est d{\'e}finie par l'{\'e}quation implicite
$$G({\wt F_c}(y),c,y) = 0~.$$
Le th{\'e}or{\`e}me des fonctions implicites est alors facile {\`a} v{\'e}rifier car
la condition qu'il requiert est exactement celle qui assure
l'existence du barycentre.

La fonction $G$ {\'e}tant $C^1$ en $(x,c,y)$ le lemme est prouv{\'e}. En
fait $F$ est, pour les m{\^e}mes raisons que pr{\'e}c{\'e}demment, $\ci$ en $c$.
\end{preuve}
\vskip5pt
{\it Preuve du cas d'{\'e}galit{\'e} {\it ii)\/}\pointir}  du th{\'e}or{\`e}me \ref{theo-princ}|

La preuve commence comme dans le paragraphe 7 de \cite{BCG1}. Pour fixer le
facteur d'homoth{\'e}tie supposons que $g$ est normalis{\'e}e de sorte que
$$\Ent(Y,g) = \Ent(\wt X,g_0) = E_0~.$$
On suppose  donc que $\vol(\rho ) = \vol (Y,g)$. Le travail porte sur
l'{\'e}tude des formes quadratiques, d{\'e}j{\`a} introduites au paragraphe pr\'ec\'edent,
\begin{align*}
h_{y,c}(\cdot,\cdot) &= \int_{\pa_F\wt X}\Big(d\ol B_{0|({\wt F_c}(y),\theta
)} (\cdot)\Big)^2 \Phi (y,c,\theta )\ d\theta \\
k_{y,c}(\cdot,\cdot) &= \int_{\pa_F\wt X} Dd\ol B_{0|({\wt F_c}(y),\theta
)} (\cdot,\cdot) \Phi (y,c,\theta )\ d\theta
\end{align*}
et des endomorphismes sym{\'e}triques et d{\'e}finis positifs correspondants,
$H_{y,c}$ et $K_{y,c}$ (ici, $c$ joue le r{\^o}le de $\Ent(Y,g)
(1+\varepsilon ) = E_0 (1+\varepsilon )$). La plus grande valeur
propre de $H_{y,c}$ est not{\'e}e $\mu ^c_n(y)$ et v{\'e}rifie,
$$0< \mu ^c_n(y)  < 1~,$$
en effet, l'endomorphisme sym{\'e}trique $H_{y,c}$ est de $\bar g_0$-trace
{\'e}gale {\`a} $1$ et est d{\'e}fini positif. On rappelle {\'e}galement que  $\tr  (K_{y,c})
= \Ent (\wt X,g_0) = E_0$ (ceci car $\Phi $ est normalis{\'e}e).

\bigskip
{\bf 1{\`e}re {\'e}tape~: convergence presque s{\^u}re de $H_{y,c}$.}

La preuve de l'in{\'e}galit{\'e} {\it i)\/} du th{\'e}or{\`e}me \ref{theo-princ} consiste (comme
dans le paragraphe 3) {\`a} montrer que, 
$$\forall y \in \wt Y, ~~\forall
c>E_0, ~~|\Jac {\wt F_c}(y)| \le \Big({c\over E_0}\Big)^n~.$$

\begin{lemme}
Il existe une suite $c_k$ tendant vers $E_0$, telle que
$\Jac {\wt F_{c_k}}(y) \build \la_{k\to + \infty}\fin 1$ presque s{\^u}rement
sur $\wt Y$.
\end{lemme}

\begin{preuve} Comme dans \cite{BCG1}, lemme 7.3, posons $f_c(y)
= \Jac {\wt F_c}(y) -1$ et $f^{\pm}_c = \sup (0,\pm f_c)$~; la fonction
$f^+_c$ tend uniform{\'e}ment vers $0$ lorsque $c$ tend vers $E_0$ car,
$$\forall y\in \wt Y,~~ 0\le f^+_c(y) \le \Big({c\over E_0}\Big)^n
-1~.$$
Par ailleurs, pour tout $c> E_0$,
\begin{align*}
\vol(\rho ) = \int_Y {\wt F_c}^*(\omega _0) 
&= \int_Y \Jac {\wt F_c}(y)\ dv_g\\
&\le \Big({c\over E_0}\Big)^n \vol (Y,g) - \int_Y f^-_c\ dv_g
\end{align*}
l'hypoth{\`e}ses $\vol(\rho ) = \vol(Y,g)$ implique que $f^-_c$ tend vers
$0$ dans $L^1(Y,g)$ lorsque $c$ tend vers $E_0$, d'o{\`u}  l'existence d'une
sous-suite $c_k$ telle que $f^-_{c_k}$ tende vers z{\'e}ro presque
s{\^u}rement.
\end{preuve}

Lorsque $(\wt Y,\tilde g) = (\wt X, \bar g_0)$ et la mesure $\mu _0$ est la
mesure canonique du bord de Furstenberg $\pa_F \wt X$, alors
l'endomorphisme  $H_{y,\varepsilon }$ prend une forme particuli{\`e}re~; en effet,
pour tout $x \in \wt X$ et pour $\varepsilon  =0$
$$
H_{x,0} =\begin{pmatrix}
{1\over pn_1} I_1 &0&0&0\\
0 &{1\over pn_2}I_2 &0&0\\
0 &0&\ddots &0\\
0&0&0&{1\over pn_p}I_p\\
\end{pmatrix},
$$
o{\`u} $x = (x_1\ld x_p)$ et $I_k$ d{\'e}signe l'identit{\'e} de $T_{x_k}\wt
X_k$. D{\'e}sormais nous noterons $H_0$ l'endomorphisme $H_{x,0}$.
De m{\^e}me, les termes $K_i$ (voir le paragraphe pr{\'e}c{\'e}dent) qui se calculent
en fonctions de $H_i = {1\over pn_i} I_i$ et valent $K_i = {E_i\over
\sqrt p n_i} I_i$.
Nous noterons $K_0$ l'endomorphisme  $K_{x,0}$ correspondant.

\`A partir de maintenant nous consid\`ererons une suite $c_k\build \la_{k\to +\infty}\fin E_0$ telle que $\Jac {\wt F_{c_k}}(y) \build \la_{k\to + \infty}\fin 1$  presque s\^urement en $y\in\wt Y$.
\begin{lemme}\label{4.6}
Pour presque tout $y \in \wt Y$, $\lim_{k\to +\infty} H_{y,c_k} = H_0$.
\end{lemme}

\begin{preuve} Pour tout $y \in \wt Y$ et pour tout $x>
E_0$
$$|\Jac {\wt F_c}(y)| \le \Big({c\over \sqrt n}\Big)^n {(\det
H_{y,c})^{1/2}\over \det (K_{y,c})} \le \Big({c\over E_0}\Big)^n~.$$
Soit $y \in Y$ tel que $|\Jac {\wt F_{c_k}}(y)| \build \la_{k\to +\infty}\fin 1$, la
quantit{\'e} ${(\det H_{y,c_k})^{1/2}\over \det( K_{y,c_k})}$ tend vers sa
valeur maximale, {\`a} savoir $\big({\sqrt n\over E_0}\big)^n$. On
rappelle que $\prodd^p_{i=1} a^{n_i}_i = 1$ (voir le paragraphe 1).

Par une preuve en tout point analogue \`a celle donn\'ee dans l'appendice B, proposition
B5 de \cite{BCG1}, nous montrons l'existence d'une constante $A>0$
telle que
$${(\det H_{y,c})^{1/2}\over \det( K_{y,c})} \le \Big({\sqrt n\over
E_0}\Big)^n \big( 1 - A \|H_{y,c}-H_0\|^2_{\bar g_0}\big)$$
de sorte que
$$\|H_{y,c}-H_0\|^2 _{\bar g_0} \le {1\over A} \Big( 1 -
\Big({E_0 \over c}\Big)^n |\Jac {\wt F_c}(y)|\Big)$$
et, si $|\Jac {\wt F_{c_k}}(y)| \build \la_{k\to +\infty}\fin 1$, alors
$$H_{y,c_k} \build \la_{k\to +\infty}  \fin H_0~.$$
\end{preuve}

\bigskip
{\bf 2{\`e}me {\'e}tape~: convergence uniforme de $H_{y,c_k}$ vers $H_0$.}

Nous reprenons les {\'e}tapes de la preuve du cas d'{\'e}galit{\'e} de \cite{BCG1},
paragraphe 7.

Soit $c_k \build\la_{k\to +\infty}\fin E_0$ une sous-suite telle que
$\Jac {\wt F_{c_k}} \la 1$ presque s{\^u}rement et $H_{y,c_k}$ tende presque
s{\^u}rement vers $H_0$. Pour simplifier les notations nous
utiliserons l'indice $k$ en lieu et place de l'indice $c_k$.

\begin{lemme}\label{4.7}
Soient $y$ et $y'$ deux points de $\wt Y$ tels que $\mu
^k_n \le 1-{1\over n}$ en tout point d'une $g$-g{\'e}od{\'e}sique minimisante
$\alpha $ qui joint $y$ {\`a} $y'$, alors
$$d_{\bar g_0}({\wt F_k}(y), {\wt F_k}(y')) \le K_1 d_g(y,y')~.$$
On rappelle que $\mu^k_n(y)=\mu^{c_k}_n(y)$ est la plus grande valeur propre de $H_{y,c_k}$. 
\end{lemme}

\begin{preuve} On tire, comme dans le paragraphe 3, de
l'{\'e}quation implicite qui d{\'e}finit ${\wt F_k}$, pour tous $u \in T_y \wt Y$
et $v\in T_{{\wt F_k}(y)} \wt X$,
\begin{align*}
\bar g_0\big(K_{y,k} D_y {\wt F_k}(u),v\big) 
&= \int_{\pa_F \wt X} d\ol B_0{}_{|
({\wt F_k}(y),\theta )} (v) \ d\Phi _k{}_{|(y,\theta )}(u)\ d\theta\\
&= 2 \int_{\pa_F \wt X} d\ol B_0{}_{| ({\wt F_k}(y),\theta )} (v)
\sqrt{\Phi _k} (y,\theta ) \ d\sqrt {\Phi _k}_{| (y,\theta
)}(u)\ d\theta \\
&\le 2 \bar g_0\big(H_{y,k}(v),v\big)^{1/2} \Big( \int_{\pa _F\wt X}
\big(d\sqrt{\Phi _k}{}_{|(y,\theta )}(u)\big)^2 d\theta
\Big)^{1/2}~.\end{align*}
Un calcul imm{\'e}diat montre que
$$\Big(\int_{\pa_F \wt X} \big(d\sqrt {\Phi _k}_{| (y,\theta )}
(u)\big)^2\ d\theta \Big)^{1/2} \le {c_k\over 2} g(u,u)^{1/2}~.\qquad (\ast )$$
Si $u$ et $v$ sont de norme 1, dans leur espace respectif, alors
\begin{align*}
\bar g_0\big(K_{y,k}  D_y {\wt F_k}(u),v\big) &\le c_k \bar g_0
\big(H_{y,k}(v),v\big)^{1/2}\\
&\le c_k \sqrt{\mu ^k_n(y)}~.\end{align*}

Maintenant, si $\wt X$ est un produit d'espaces sym{\'e}triques de rang 1, de
courbure comprise entre $-1$ et $-4$, il est
facile de v{\'e}rifier (voir \cite{BCG1}, appendice B) que, au sens des formes
quadratiques, pour tout $i=1,2\ld p$,
$$K_i \ge I_i - H_i \ge (1-\mu ^k_n(y)) I_i~.$$
On rappelle que $H_i$ (resp. $K_i$) d{\'e}signe la restriction de
$H_{y,k}$ (resp. $K_{y,k})$ {\`a} $\wt X_i$.
En prenant $v = {D_y {\wt F_k}(u)\over \|D_y {\wt F_k}(u)\|_{\bar g_0}}$ si $D_y{\wt F_k}(u)
\ne 0$, il vient
$$\|D_y {\wt F_k}(u)\|_{\bar g_0}\le c_k {{\sqrt \mu ^k_n(y)}\over 1- \mu
^k_n(y)}\qquad (\ast \ast )$$
(si $D_y{\wt F_k} (u) = 0$, l'in{\'e}galit{\'e} est trivialement vraie). Soit $\alpha
$ la $g$-g{\'e}od{\'e}sique de $y$ {\`a} $y'$ le long de laquelle $\mu ^k_n(\alpha
(t)) \le 1 - {1\over n}$, on a, pour tout $u \in T_{\alpha (t)} \wt Y$,
de norme 1
$$\|D_{\alpha (t)} {\wt F_k}(u)\|_{\bar g_0} \le 2n E_0 = K_1$$
(si $k$ est assez grand pour que $c_k \le 2E_0$). Par le th{\'e}or{\`e}me des
accroissements finis
$$d_{\bar g_0} \big( {\wt F_k}(y),{\wt F_k}(y')\big) \le K_1 d_g(y,y')~.$$
\end{preuve}

\begin{lemme}\label{4.8}
Avec les m\^emes notations que pr\'ec\'edemment, si $P$ d{\'e}signe le transport
parall{\`e}le de ${\wt F_k}(y)$ {\`a} ${\wt F_k}(y')$ le long de la $\bar g_0$-g{\'e}od{\'e}sique
minimisante qui les joint, on a
$$\|h_{y,k}\circ P - h_{y',k}\|_{\bar g_0} \le K_2\big[ d_g(y,y') +
d_{\bar g_0}
({\wt F_k}(y),{\wt F_k}(y'))\big]~.$$
\end{lemme}

\begin{preuve} Nous d{\'e}signons par $\beta (t)$ l'unique
$\bar g_0$-g{\'e}od{\'e}sique, qui est minimisante, allant de ${\wt F_k}(y)$ {\`a} ${\wt F_k}(y')$
et par $Z$ un champ de vecteurs parall{\`e}le, le long de $\beta $, de
norme 1. Pour simplifier, posons $Z_1 = Z({\wt F_k}(y))$ et $Z_2 =
Z({\wt F_k}(y'))$. Alors
\begin{align*}
h_{y',k}(Z_2,Z_2) - h_{y,k}&(Z_1,Z_1)\\
&= \int_{\pa_F\wt X} \Big(d\ol B_0{}_{| ({\wt F_k}(y'),\theta )}(Z_2)\Big)^2
\Phi _k(y',\theta )\ d\theta \\
&\kern1.8cm 
- \int_{\pa_F\wt X} \Big(d\ol B_0{}_{| ({\wt F_k}(y),\theta )}(Z_1)\Big)^2
\Phi _k(y,\theta )\ d\theta \\
&=  \int_{\pa_F\wt X} \Big[\Big(d\ol B_0{}_{| ({\wt F_k}(y'),\theta )}(Z_2)\Big)^2-
\Big(d\ol B_0{}_{| ({\wt F_k}(y),\theta )}(Z_1)\Big)^2\Big]\Phi _k(y',\theta )\ d\theta \\
&\kern1.8cm  
+ \int_{\pa_F\wt X} \Big(d\ol B_0{}_{| ({\wt F_k}(y),\theta )}(Z_1)\Big)^2
\big(\Phi  _k(y',\theta )- \Phi  _k(y,\theta )\big)\ d\theta
~.\\
\end{align*}

Des formules explicites de $Dd\ol B_0$ et du fait que
$\|d\ol B_0{}_{|(x,\theta )} (\cdot)\|_{\bar g_0} \le 1$, nous tirons
l'in{\'e}galit{\'e}
$$|\Big(d\ol B_0{}_{| ({\wt F_k}(y'),\theta )}(Z_2)\Big)^2-
\Big(d\ol B_0{}_{| ({\wt F_k}(y),\theta )}(Z_1)\Big)^2| \le K'_2
d_{\bar g_0}\big({\wt F_k}(y'),{\wt F_k}(y)\big)~.$$
De m{\^e}me, comme $\Phi _k(y,\cdot)$ est de norme 1 dans $L^1(\pa_F\wt
X,d\theta )$ et en utilisant l'in\'equation $(\ast )$
\begin{align*} \int_{\pa_F\wt X} \Big(d\ol B_0&{}_{| ({\wt F_k}(y),\theta )}(Z_1)\Big)^2
\Big[ \Big( \sqrt{\Phi  _k(y',\theta )}\Big)^2 - \Big(\sqrt{\Phi
_k(y,\theta )}\Big)^2\Big]\ d\theta\\
&\le \Big(\|\sqrt{\Phi _k}(y,\cdot) - \sqrt{\Phi _k}(y',\cdot)\|_{L^2(\pa_F\wt
X)}\Big)  \Big(\|\sqrt{\Phi _k}(y,\cdot) + \sqrt{\Phi _k}(y',\cdot)\|_{L^2(\pa_F\wt
X)}\Big)\\
&\le c_k d_g(y,y')~.\end{align*}
Le lemme d{\'e}coule de l'addition de ces in{\'e}galit{\'e}s.
\end{preuve}

\begin{lemme}
La suite $H_{k,y}$ converge uniform{\'e}ment par rapport \`a $y\in \wt
Y$ vers $H_0$ lorsque $k$ tend vers $+\infty$.
\end{lemme}

\begin{preuve} Le comportement de $H_k$ vis-{\`a}-vis de
l'action de $\Gamma $ sur $\wt Y$ montre qu'il suffit de prouver la
convergence uniforme sur un domaine fondamental $D\subset \wt Y$. Le
th{\'e}or{\`e}me d'Egoroff (\cite{Fed}, p.~77) et le lemme \ref{4.6} attestent que, pour
tout $\eta >0$, il existe un ensemble mesurable $K$ tel que
\begin{itemize}
\item[i)] $\vol_g(D\setminus K) < \eta $~;

\item[ii)] sur $K$, la suite $y\mapsto H_{k,y}$ converge uniform{\'e}ment
vers $H_0$.
\end{itemize}
\noindent Fixons $\varepsilon >0$ petit, on peut choisir $\eta $ tel
que $D\setminus K$ ne contienne aucune $g$-boule de rayon $\varepsilon $,
car, en effet, le volume d'une telle boule sur $\wt Y$ est minor{\'e} (la
m{\'e}trique de $\wt Y$ est p{\'e}riodique). On choisit aussi $N\in \textbf{N}$ de
sorte que
\begin{itemize}
\item[i)] pour tout $k\ge N$, $E_0 <c_k< E_0+\varepsilon $~;

\item[ii)] pour tout $k\ge N$ et pour tout $y\in K$,
$\|H_{y,k} - H_0\|_{\bar g_0} < \varepsilon $.
\end{itemize}
\noindent Par ailleurs, si $y \notin K$, $d_g(y,K) < \varepsilon
$. Rappelons que les valeurs propres de $H_0$ sont les nombres
${1\over pn_i}$, $i=1,2\ld p$. Posons $K_3 = K_2(K_1+1) +1$ et supposons $\varepsilon $ assez
petit pour que $K_3 \varepsilon  \le 1 - \dis\sup_i\big({1\over
pn_i}\big)- {1 \over n}$. Nous allons
montrer que si $k\ge N$, alors
$$\forall y \in D,~~ \|H_{y,k} - H_0\|_{\bar g_0} < K_3 \varepsilon ~.$$
Si ce n'est pas vrai, il existe $y' \in D$ tel que
$$\|H_{y',k}- H_0\|_{\bar g_0} \ge K_3 \varepsilon ~,$$
soit alors $y \in K$ tel que $d_g(y',y) < \varepsilon $. Par
continuit{\'e} de l'application $y \mapsto H_{y,k}$, il existe un premier
point $y''$ sur le segment g{\'e}od{\'e}sique $[y,y']$ tel que 
$\|H_{y'',k} - H_0\| = K_3 \varepsilon$. Le choix 
de $K_3$ montre que, sur le segment g{\'e}od{\'e}sique $[y,y'']$,
$$\mu ^k_n \le \sup_i \Big({1\over pn_i}\Big) + K_3\varepsilon  \le  1-{1\over n}~.$$

D'apr{\`e}s les lemmes \ref{4.7} et \ref{4.8} ceci conduit {\`a}
$$\|h_{y,k}\circ P - h_{y'',k}\|_{\bar g_0} \le K_2(K_1+1)\varepsilon $$
et comme $\|H_{y,k} - H_0\|_{\bar g_0} < \varepsilon $ 
ceci conduit {\`a}
$$\|H_{y'',k} - H_0\|_{\bar g_0} < \big(K_2(K_1 + 1)+1\big)
\varepsilon  = K_3 \varepsilon $$
qui est une contradiction.
\end{preuve}

Remarquons que la convergence uniforme de $H_{y,k}$ vers $H_0$ implique la convergence uniforme de $K_{y,k}$ vers $K_0 $.

\bigskip
{\bf 3{\`e}me {\'e}tape~: convergence uniforme d'une sous-suite de ${\wt F_k}$.}

\begin{lemme}
Il existe une sous-suite de la suite ${\wt F_k}$ qui converge
uniform{\'e}ment vers une application ${\wt F} : \wt Y \to \wt X$ continue et
{\'e}quivariante.
\end{lemme}

\begin{preuve} Pour $\varepsilon >0$ donn{\'e}, il existe
$M\in \textbf{N}$ tel que si $k\ge M$
$$\forall y \in \wt Y,~~ \|H_{y,k} - H_0\|_{\bar g_0} < \varepsilon ~.$$
D'o{\`u}
$$H_{y,k} \le H_0+ \varepsilon  I$$
et par une remarque pr{\'e}c{\'e}dente
$$K_{y,k} \ge K_0 - \varepsilon  I~.$$
Ces deux in{\'e}galit{\'e}s {\'e}tant {\`a} comprendre au sens des formes
quadratiques. On d{\'e}duit alors, avec $(\ast \ast )$, qu'il existe un nombre r{\'e}el $C>0$ tel
que,  pour tout $y
\in \wt Y$ et $u \in T_y \wt Y$,
$$\|D_y {\wt F_k}(u)\|_{\bar g_0} \le  C$$
(si $\varepsilon $ est assez petit).

La suite d'application ${\wt F_k} : \wt Y \to \wt X$ est donc {\'e}quicontinue.

Supposons qu'il existe $y_0$ tel que ${\wt F_k}(y_0)$ ne reste dans aucun
compact. Quitte {\`a} extraire une sous-suite, on peut supposer que
${\wt F_k}(y_0) \build \la_{k\to +\infty}\fin \theta \in \pa \wt X$ (le bord
g{\'e}om{\'e}trique de $\wt X$). Pour tout $y \in \wt Y$, alors
$$d_{\bar g_0} ({\wt F_k}(y),{\wt F_k}(y_0)) \le C\ d_g(y,y_0)$$
de sorte que ${\wt F_k}(y) \build \la_{k\to +\infty}\fin \theta $ par
d{\'e}finition du bord g{\'e}om{\'e}trique de $\wt X$. L'{\'e}quivariance de ${\wt F_k}$
donne
$${\wt F_k}(\gamma y_0) = \rho (\gamma ) {\wt F_k} (y_0)$$
et donc en passant {\`a} la limite en $k$
$$\theta  = \rho (\gamma )\theta $$
c'est-{\`a}-dire, la repr{\'e}sentation $\rho $ fixe un point de $\pa\wt X$.

\begin{lemme}
\label{pointfixe}
Si $\rho $ fixe un point $\theta _0$ de $\pa \wt X$,
alors $\vol(\rho ) =0$.
\end{lemme}

\begin{preuve} Soit $\ol B_0(\cdot,\theta _0)$ la fonction de
Busemann d{\'e}finie par le point $\theta _0 \in \pa \wt X$. Supposons d'abord que $\theta_0$ est dans le bord de F{\"u}rstenberg. Posons
$$Z(x) = \nabla \ol B_0(x,\theta _0)$$
alors le champ de vecteurs $Z$ est invariant par $\rho $. En effet,
l'{\'e}galit{\'e}
$$\ol B_0 \big(\alpha (x),\theta _0\big) = \ol B_0\big( x,\alpha ^{-1}(\theta
_0)\big) + \ol B_0\big(\alpha (O _0),\theta _0\big)$$
pour $\alpha  \in \Isom (\wt X)$, conduit {\`a}
$$\ol B_0\big(\rho (\gamma )(x),\theta _0\big) = \ol B_0(x,\theta _0) +
\ol B_0\big(\rho (\gamma )(O_0),\theta _0\big)$$
pour tout $\gamma  \in \Gamma $~; ce qui donne en diff{\'e}renciant
$$Z\big(\rho (\gamma )(x)\big) = \rho (\gamma )(Z(x))~.$$
Par ailleurs, pour tout $x \in \wt X$
$$\divv(Z)(x) = \Delta \big( \ol B_0(\cdot,\theta )\big) = E_0~.$$
Donc la forme diff{\'e}rentielle $\omega  =\divv(Z)\omega _0 =  E_0 \omega _0$
est invariante par $\rho (\gamma )$, pour tout $\gamma \in \Gamma
$. En cons{\'e}quence, pour $c>E_0$, ${\wt F}^*_c(\omega )$ est invariante par $\gamma $,
pour tout $\gamma  \in \Gamma $. La d{\'e}finition de la divergence
conduit {\`a} l'{\'e}galit{\'e}
$$\divv(Z) \omega _0 = -d \big(i(Z)\cdot \omega _0\big)$$
o{\`u} $i(Z)\cdot \omega _0$ d{\'e}signe le produit int{\'e}rieur de $\omega _0$
par le champ de vecteurs $Z$. D'o{\`u}
\begin{align*}
{\wt F}^*_c (\omega) &= -{\wt F}^*_c\big(d(i(Z)\cdot \omega _0)\big)\\
&= -d\big({\wt F}^*_c(i(Z)\cdot \omega _0)\big)\\\end{align*}
et
$$\vol(\rho ) = \int_Y {\wt F}^*_c (\omega _0) = {1\over E_0} \int_Y {\wt F}^*_c
(\omega ) =0~.$$
Si $\theta_0$ n'est pas dans le bord de F{\"u}rstenberg la m\^eme preuve peut \^etre faite car 
$$\divv(Z)(x) = \Delta \big( \ol B_0(\cdot,\theta )\big) \ne 0~.$$
\end{preuve}

Puisque nous sommes dans le cas d'{\'e}galit{\'e}, $\vol(\rho ) \ne 0$, et la
suite ${\wt F_k}(y_0)$ reste donc dans un compact de $\wt X$. On peut alors
appliquer le th{\'e}or{\`e}me d'Ascoli pour d{\'e}duire qu'il existe une
sous-suite, not{\'e}e encore ${\wt F_k}$, qui converge uniform{\'e}ment sur
$D\subset \wt Y$ vers une application continue $F : D \to \wt
X$. L'{\'e}quivariance de ${\wt F_k}$, pour tout $k$, montre que ${\wt F_k}$ converge
uniform{\'e}ment sur $\wt Y$ et que la limite $F$ est {\'e}galement
{\'e}quivariante.

\bigskip
{\bf 4{\`e}me {\'e}tape~: $\wt F$ est une isom{\'e}trie.}
\vskip-10pt

\begin{lemme}
L'application ${\wt F} : (\wt Y,\tilde g)  \to (\wt X,g_0)$ contracte les
distances, c'est-{\`a}-dire, pour tout $y$, $y'$ dans $Y'$
$$d_{g_0} \big({\wt F}(y),{\wt F}(y')\big) \le d_{g} (y,y')$$
et $D{\wt F}(y)$ est une isom{\'e}trie entre $(T_y \wt Y,g)$ et $(T_{{\wt F}(y)} \wt X,g_0)$
pour presque tout $y \in \wt Y$.
\end{lemme}

\begin{preuve} Pour $\varepsilon  >0$ donn{\'e}, on peut
choisir $k$ assez grand pour que, pour tout $y \in \wt Y$, 
\begin{align*}H_{y,k} &\le H_0+ \varepsilon I\\
K_{y,k} &\ge K_0 - \varepsilon I ~.\\\end{align*}
Alors, l'in{\'e}galit{\'e} 2.5 nous conduit {\`a} l'estimation suivante, pour
$u\in T_y \wt Y$ et $v \in T_{{\wt F_k}(y)} \wt X$~:
$$\bar g_0\big(K_{y,c}(D{\wt F_k}(y)(u)),v\big) \le (1+\varepsilon ) E_0
\big(\bar g_0 (H_{y,c}(v),v)\big)^{1/2} \big(h'_{y,k}
(u,u)\big)^{1/2}~.$$
On rappelle que la d{\'e}finition des fonctions ${\wt F_k}$ est ind{\'e}pendante des
coefficients choisis pour d{\'e}finir la m{\'e}trique de r{\'e}f{\'e}rence, 
\cad
qu'elle donne la m{\^e}me fonction qu'on utilise $\bar g_0 = \opp^p_{i=1}
g^i_0$ ou bien $g_0 = \opp^p_{i=1}a^2_i g^i_0$. Nous avons choisi
d'utiliser $\bar g_0$ pour d{\'e}finir le barycentre sur $\wt X$, en
cons{\'e}quence les matrices $H_{y,k}$, $H_0$, $K_{y,k}$ et $K_0$ sont
d{\'e}finies {\'e}galement gr{\^a}ce {\`a} la m{\'e}trique $\bar g_0$.

On rappelle {\'e}galement que $G_{i,k}$, $i = 1,2\ld p$, d{\'e}signe la
$i$-i{\`e}me composante de ${\wt F_k}$ dans la d{\'e}composition $\wt X = \wt X_1
\times \cdots \times \wt X_p$ et que $H_i$ (resp. $K_i$) d{\'e}signe la
restriction de $H_{y,k}$ (resp. $K_{y,k}$) {\`a} $T_{{\wt F}_{i,k}(y)} \wt X_i$
(ici on omet volontairement les indices $y$ et $k$ dans $H_i$ et $K_i$
afin d'all{\'e}ger les notations). Si $v = (v_1 \ld v_n)$ est tangent {\`a}
$\wt X_i$, \cad si $v_j=0$ pour tout $j \ne i$, alors, gr{\^a}ce {\`a} la
forme diagonale par blocs de $K_{y,c}$ nous obtenons
$$\bar g_0\big(K_i(D{\wt F}_{i,k}(y)(u)),v_i\big) \le (1+\varepsilon ) E_0
\big(\bar g_0 (H_iv_i,v_i)\big)^{1/2} \big(h'_{y,k}
(u,u)\big)^{1/2}$$
(on identifie, par abus de langage $v$ {\`a} sa composante $v_i$).

En utilisant l'in{\'e}galit{\'e} pr{\'e}c{\'e}dente sur $H_{y,k}$,
$$\bar g^i_0\big(K_i(D{\wt F}_{i,k}(y)(u)),v_i\big) \le (1+\varepsilon ) E_0
\Big({1\over pn_i} + \varepsilon \Big)^{1/2} \|v_i\|_{g^i_0} \big(h'_{y,k}(u,u)\big)^{1/2}~.$$

En prenant le supremum en $v_i$ de norme 1, nous obtenons,
$$\|K_{y,c}\big(D{\wt F}_{i,k} (y)(u)\big)\| _{\bar g_0} = \|K_i\big(D{\wt F}_{i,k}
(y)(u)\big)\|_{g^i_0} \le {E_0 \over \sqrt p \sqrt{n_i}}
\big(h'_{y,k}(u,u)\big)^{1/2} (1+o(\varepsilon ))~.$$
Les in{\'e}galit{\'e}s pr{\'e}c{\'e}dentes donnent encore,
$${E_i \over \sqrt{p n_i}} (1 + o(\varepsilon )) \|D{\wt F}_{i,k}
(y)(u)\|_{\bar g_0} \le {E_0\over \sqrt p \sqrt {n_i}}
\big(h'_{y,k}(u,u)\big)^{1/2} (1+o(\varepsilon ))~,$$
or les coefficients $a_i$ apparaissant dans la d{\'e}finition de la
m{\'e}trique $g_0$ valent~:
$$a_i = {E_i \sqrt n\over \sqrt{n_i} E_0}$$
d'o{\`u}, pour tout $u \in T_y \wt Y$
$$a_i \|D{\wt F}_{i,k}(y)(u)\|_{g^i_0} \le \sqrt n
\big(h'_{y,k}(u,u)\big)^{1/2} (1+o(\varepsilon ))$$
et, pour tout $u \in T_y\wt Y$
$${\wt F}^*_k g_0 (u,u) = \|D {\wt F_k}(y) (u)\|^2_{g_0} = \summ^p_{i=1}
a^2_i\|D{\wt F}_{i,k}  (y) (u)\|^2_{g^i_0} \le n h'_{y,k} (u,u) (1+o(\varepsilon ))~.$$

On peut alors calculer la trace du tenseur sym{\'e}trique ${\wt F}^*_k g_0$ par
rapport {\`a} la m{\'e}trique $g$ sur $\wt Y$ en $y \in \wt Y$.
$$\tr_g({\wt F}^*_k g_0) \le n(1+o(\varepsilon ))~.$$
En effet, on rappelle que $\tr_g(h'_{y,k}) = 1$ (voir le paragraphe
2).

Par ailleurs le d{\'e}terminant de ${\wt F}^*_k g_0$ relativement {\`a} $g$, c'est-{\`a}-dire
$|\Jac {\wt F_k}|^2$, tend presque s{\^u}rement vers $1$ sur $\wt Y$. Alors si
$A_{k,y}$ d{\'e}signe la matrice de ${\wt F}^*_k g_0$ dans une base
$g$-orthonorm{\'e}e, nous avons, pour $k$ assez grand,
$$1-\varepsilon  \le (\det A_{k,y})^{1/n} \le {1\over n} \tr (A_{k,y})
\le 1+\varepsilon $$
ce qui implique que
$$\|A_{k,y} - (\det A_{k,y})^{1/n}Id\| = o(\varepsilon )~.$$
En conclusion,  $D{\wt F_k}$ converge presque s{\^u}rement sur $\wt Y$ vers une
isom{\'e}trie.

Alors, l'application $\wt F$ est limite uniforme d'une suite
d'applications lipschitzienne  ${\wt F_k}$ dont les diff{\'e}rentielles $D{\wt F_k}$
sont uniform{\'e}ment born{\'e}es et convergent presque s{\^u}rement vers une isom{\'e}trie~; le
lemme 7.8 de \cite{BCG1} montre que, dans ce cas, l'application $F$ est
$1$-lipschitzienne. Nous ne reproduisons pas la preuve de ce fait.

L'application $\wt F$ est presque partout diff{\'e}rentiable par le th{\'e}or{\`e}me
de Rademacher et, comme elle est $1$-lipschitzienne, on a, pour
presque tout $y\in \wt Y$
$$|\Jac {\wt F}(y)|\le 1~.$$
Par ailleurs,
$$\vol(\rho ) = \int_Y {\wt F}^* (\omega _0) = \int_Y \big(\Jac {\wt F}(y)\big)
dv_g = \vol (Y,g)~.$$
D'o{\`u}, pour presque tout $y \in \wt Y$, $\Jac {\wt F}(y) =1$.

Enfin, pour presque tout $y\in \wt Y$, pour tout $u \in T_y \wt
Y$, le caract{\`e}re $1$-lipschitzien de $\wt F$ implique que
$$\|D_y {\wt F}(u)\|_{g_0} \le \|u\|_g~.$$
Ceci, combin{\'e} au fait que pour presque tout $y\in \wt Y$, $\Jac {\wt F}(y) =
1$, montre que la diff{\'e}rentielle de $\wt F$, $D_y {\wt F}$, est presque partout
sur $\wt Y$ une isom{\'e}trie (entre $T_y \wt Y$ et $T_{{\wt F}(y)} \wt X$).

Le lemme est prouv{\'e}.
\end{preuve}

\begin{lemme}
L'application $\wt F$ minimise la fonctionnelle $E_p (h) =
{1\over \vol (Y)} \int_Y\! \|Dh\|_{g,g_0}^p$ parmi toutes les applications $h$ de
$\wt Y$ dans $\wt X$, $\rho $-{\'e}quivariantes et lipschitziennes, pour
tout $p\ge n$. Ici $\|Dh\|_{g,g_0}^p$ est calcul\'ee \`a l'aide de la m\'etrique $g$ sur $\wt Y$ et $g_0$ sur $\wt X$.
\end{lemme}

\begin{preuve} Notons que, par l'\'equivariance de $h$, l'int\'egrand dans l'expression de $E_p(h)$ est invariant par $\Gamma$ et est donc une quantit\'e d\'efinie sur $Y$. Si $\{e_i\}$ est une base $g$-orthonorm{\'e}e en $y\in \wt
Y$.
$$\|Dh(y)\|_{g,g_0} = \Big( {1\over n} \summ^n_{i=1} \|Dh(y)
(e_i)\|^2_{g_0}\Big)^{1/2}~.$$
Cette quantit{\'e} est d{\'e}finie pour presque tout $y\in \wt Y$. On a donc,
pour presque tout $y\in \wt Y$,
$$|\Jac h(y)|^{p/n} \le \|Dh(y)\|^p_{g,g_0}$$
pour tout $p\ge 0$. Maintenant si $p\ge n$
\begin{align*}
1 = \Big( {\vol (\rho )\over \vol (Y)}\Big)^{p/n} = \Big( {1\over
\vol (Y)} \int_Y \Jac h(y) dv_g(y)\Big)^{p/n} &\le {1\over \vol(Y)}
\int |\Jac h(y)|^{p/n} dv_g(y)\\
& \le E_p (h)~.\end{align*}

Si $h$ est remplac{\'e}e par $\wt F$, en utilisant le fait que $D{\wt F}(y)$ est une
isom{\'e}trie pour presque tout $y \in \wt Y$, il vient
$$1 = E_p({\wt F}) \le E_p(h)~.$$
\end{preuve}

\begin{corollaire}
L'application ${\wt F}$ est de classe $\ci$.
\end{corollaire}

\begin{preuve} En fait, nous prouvons que ${\wt F}$ est
harmonique, la r{\'e}gularit{\'e} s'en d{\'e}duit.

De mani{\`e}re heuristique nous pouvons dire que l'{\'e}quation d'Euler
associ{\'e}e {\`a} la fonctionnelle $E_p$, $p\ge n$ s'{\'e}crit
$$\divv\big( \|D{\wt F}\|^{p-2}_{g,g_0} D{\wt F}\big) = 0$$
o{\`u} la divergence est {\`a} comprendre comme celle d'une 1-forme sur $\wt
Y$ {\`a} valeurs dans $T\wt X$ (voir \cite{Eel-Lem1}, page~6). Mais ${\wt F}$ a une
diff{\'e}rentielle qui est presque partout une isom{\'e}trie, de sorte que
$\|D{\wt F}\|_{g,g_0} =1$ presque partout sur $\wt Y$, et l'{\'e}quation devient
$$\divv(D{\wt F}) =0$$
c'est-{\`a}-dire ${\wt F}$ est harmonique.

Plus pr{\'e}cis{\'e}ment, $D{\wt F}$ est interpr{\'e}t{\'e}e comme une 1-forme sur $\wt Y$ {\`a}
valeurs dans le fibr{\'e} ${\wt F}^{-1}(T\wt X)$, c'est-{\`a}-dire un {\'e}l{\'e}ment de $C^0_\rho
\big(T^*(\wt Y) \otimes {\wt F}^{-1} (T\wt X)\big)$, qui est de plus $\rho
$-{\'e}quivariante (voir \cite{Eel-Lem1}, page~8)~; soit alors $Z$ un champ de
vecteurs $\ci$ le long de ${\wt F}$, qui satisfait {\'e}galement la relation
de $\rho $-{\'e}quivariance ad{\'e}quate, c'est-{\`a}-dire qui est un {\'e}l{\'e}ment de $\ci_\rho \big(\wt
Y,{\wt F}^{-1}(T\wt X)\big)$; alors il existe une variation {\`a} un param{\`e}tre de ${\wt F}$,
not{\'e}e ${\wt F}_t$, $\rho $-{\'e}quivariante, telle que
$$\forall y \in \wt Y,~~ {d\over dt} {\wt F}_t (y)_{|t=0} = Z(y)$$
(voir \cite{Eel-Lem2}, page~397).

Comme ${\wt F}$ minimise $E_p$, pour $p\ge n$, on a
$${d\over dt}_{|t=0} E_p({\wt F}_t) = 0$$
c'est-{\`a}-dire,
$${d\over dt}_{|t=0} {1\over \vol(Y)} \int_Y \|D{\wt F}_t (y)\|^p_{g,g_0}
dv_g(y)=0$$
mais
\begin{align*}
{d\over dt}\Big( {1\over \vol(Y)} \int_Y \|D{\wt F}_t (y)\|^p_{g,g_0} &dv_g(y)\Big)\\
&= {1\over \vol(Y)}\! \int_Y {d\over dt} \big(\|D{\wt F}_t (y)\|^p_{g,g_0}\big) dv_g(y)\\
&= {1\over \vol(Y)} \!\int_Y {p\over 2} \|D{\wt F}_t (y)\|^{p-2}_{g,g_0}{d\over dt} \big(\|D{\wt F}_t (y)\|^2_{g,g_0} \big) dv_g(y).\\\end{align*}
En $t=0$, comme $\|D{\wt F}(y)\|_{g,g_0} = 1$ pour presque tout $y \in \wt Y$, on a
$$0 = {p/2\over \vol(Y)} \int_Y {d \over dt}_{| t=0}
\big(\|D{\wt F}_t(y)\|^2_{g,g_0}\big) dv_g(y)$$
c'est-{\`a}-dire, ${\wt F}$ est un point critique de la fonctionne $E_2$. L'application
${\wt F}$ est donc faiblement harmonique (au sens des distributions, voir
\cite{Eel-Lem2}, page~397). D'apr{\`e}s les th{\'e}or{\`e}mes de r{\'e}gularit{\'e} classiques
(voir \cite{Eel-Lem2}, 3.10, page~397), ${\wt F}$ {\'e}tant continue, elle est de
classe~$\ci$.
\end{preuve}
\begin{remarque}
Nous avons montr{\'e} que ${\wt F}$ est un point critique
de $E_2$, mais en fait elle minimise cette fonctionnelle car l'espace
{\'e}tant de courbure n{\'e}gative ou nulle la fonctionnelle $E_2$ est
convexe.
\end{remarque}

Nous pouvons alors terminer la preuve du th{\'e}or{\`e}me \ref{theo-princ}~{\it
ii)}. L'application ${\wt F}$ a une diff{\'e}rentielle $D{\wt F}(y)$ qui est continue
en $y$ et est donc une isom{\'e}trie pour tout $y \in \wt Y$~; la vari{\'e}t{\'e}
$\wt Y$ {\'e}tant connexe et compl{\`e}te, $\wt X$ {\'e}tant connexe et simplement
connexe nous d{\'e}duisons de cela que ${\wt F}$ est une isom{\'e}trie surjective de
$\wt Y$ sur $\wt X$ (c'est en effet un exercice classique, voir
\cite{GHL}, 2.108, exercice a), page 97). En particulier $\rho (\Gamma )$
est un sous-groupe discret cocompact de $\Isom(\wt X)$ agissant sans
points fixes et la repr{\'e}sentation $\rho $ est injective.
\end{preuve}
\begin{remarques}
\label{3.18}
\begin{itemize}
\item[i)] Le lemme \ref{pointfixe} peut s'{\'e}tendre et donne lieu {\`a} la proposition
suivante~:
\end{itemize}
\begin{proposition}
\label{radon-inv}
S'il existe une mesure de Radon finie et non nulle $\mu $,
d{\'e}finie sur $\pa \wt X$, invariante par $\rho (\Gamma )$, alors
$\vol(\rho ) = 0$.
\end{proposition}

\begin{preuve} La preuve est identique {\`a} celle du lemme \ref{pointfixe}, en posant
$$Z(x) = \int_{\pa \wt X} \nabla B(x,\theta ) d\mu (\theta
)~.$$
\end{preuve}

\begin{itemize}
\item[ii)] Par ailleurs, si $\vol (\rho ) \ne 0$ le groupe $\rho
(\Gamma )$ ne peut pas fixer (globalement) un sous-espace strict et
totalement g{\'e}od{\'e}sique de $\wt X$, car, sinon, nous pourrions choisir
une application {\'e}quivariante $f$ {\`a} valeurs dans ce sous-espace, et la
chute de dimension entra{\^\i}nerait que $\vol(\rho ) = 0$, une
contradiction. En utilisant le crit{\`e}re g{\'e}om{\'e}trique {\'e}nonc{\'e} dans \cite{Lab},
nous montrons donc ({\`a} l'aide de la remarque {\it ii)\/} et de la
proposition \ref{radon-inv}) que
\end{itemize}

\begin{proposition}\label{reductif}
Si $\vol (\rho ) \ne 0$ alors $\rho (\Gamma )$
est r{\'e}ductif.
\end{proposition}

On rappelle que $\rho (\Gamma)$ est dit r\'eductif si son adh\'erence de Zariski l'est, c'est-\`a-dire si cette derni\`ere a un radical unipotent trivial.

Notons que dans \cite{Lab}, la r{\'e}ductivit{\'e} de $\rho (\Gamma )$ est prouv{\'e}e
{\^e}tre une condition n{\'e}cessaire et suffisante {\`a} l'existence d'une
application harmonique $\rho $-{\'e}quivariante. 

\end{remarques}
 
Enfin, le th{\'e}or{\`e}me \ref{theo-princ}
conduit au

\begin{corollaire}
\label{3.21}
Si $\rho $ est une repr{\'e}sentation de $\Gamma  =
\pi _1(Y)$ dans $\Isom (\wt X,g_0)$, o{\`u} $Y$ est une vari{\'e}t{\'e} compacte,
alors
$$\minvol (Y) \ge \Big( {\Ent(\wt X,g_0)\over n-1}\Big)^n \vol(\rho
)~.$$
\end{corollaire}

\begin{preuve} On rappelle que
$$\minvol(Y) = \inf \{\vol (Y,g) \\ g \hbox{~m{\'e}trique sur $Y$  telle que~}
|K_g| \le 1\}$$
et que si la courbure sectionnelle $K_g$ de la m{\'e}trique $g$ v{\'e}rifie
$K_g \ge -1$ alors on a $\Ent(Y,g) \le n-1$ (voir \cite{BCG1}).
\end{preuve}

\begin{remarque}

{\it i)\/}~~ En particulier, s'il existe une repr{\'e}sentation $\rho $
telle que $\vol(\rho ) \ne 0$ alors $\minvol(Y) >0$.

{\it ii)\/}~~ On pourrait remplacer le volume minimal $\minvol(Y)$ par
$$\minvol_{\ric} (Y) = \inf\big\{\vol(Y,g)\\ \ric_g \ge -
(n-1)g\big\}~.$$
\end{remarque}

\section{Applications}

Dans ce paragraphe nous nous int{\'e}ressons au cas o{\`u} $(\wt Y, \tilde g)$
est elle-m{\^e}me un produit fini d'espaces sym{\'e}triques simplement connexe
de courbure strictement n{\'e}gative. Comme pr{\'e}c{\'e}demment  un tel espace
sera not{\'e} $(\wt X, \tilde g_0)$, o{\`u} $g_0$ est la m{\'e}trique d{\'e}finie au
paragraphe 2 et qui minimise l'entropie. De m{\^e}me, $\Gamma$ d{\'e}signe un
r{\'e}seau cocompact et sans torsion de $\ixg$, et
$\rho $ est un morphisme 
$$\rho : \Gamma\la\ixg.$$
Des exemples de telles repr{\'e}sentations sont rares et le but de ce
paragraphe est, en particulier, de rappeler quelques unes des
constructions classiques.

Dans cette situation, le th{\'e}or{\`e}me \ref{theo-princ} s'{\'e}crit 
$$\vol(\rho )\le\vol(X, g_0)$$
o{\`u} $X=\wt X/\Gamma$. L'{\'e}galit{\'e}, dans cette in{\'e}galit{\'e}, n'a lieu que si
et seulement si $(\wt X/\rho(\Gamma ), g_0)$ est une vari{\'e}t{\'e}
isom{\'e}trique {\`a} $(X,g_0)$, c'est-{\`a}-dire si $\rho (\Gamma )$ est un
r{\'e}seau cocompact de $\ixg$. Nous r{\'e}pondons, dans
ce paragraphe {\`a} la question~:       

\begin{question}
Existe-t-il des repr{\'e}sentations, comme ci-dessus, telles
que $0<\vol(\rho ) < \vol(X,g_0)$~?
\end{question}

Rappelons qu'un r\'eseau $\Gamma$ dans un groupe de Lie $G$,
semi-simple connexe sans facteur compact est dit r\'eductible si $G$
poss\`ede des sous-groupes normaux $H$ et $H'$ tels que $G=H.H'$,
$H\cap H'$ est discret et $\Gamma /(\Gamma \cap H). (\Gamma \cap
H')$ est fini (voir \cite{Rag} page 86). $\Gamma$ est dit
irr\'eductible s'il n'est pas r\'eductible

Alors, lorsque $\Gamma $ est irr{\'e}ductible, le th{\'e}or{\`e}me de
super-rigidit{\'e} de Margulis (\cite{Mar}, chapitre VII) fournit une
r{\'e}ponse n\'egative compl{\`e}te \`a la question ci-dessus.

\begin{proposition}
\label{irreductible}
Avec les notations ci-dessus, si $\Gamma $ est
irr{\'e}ductible et $\vol(\rho ) \ne 0$ alors $\rho (\Gamma )$ est un
r{\'e}seau cocompact de $\ixg$ et donc $\vol(\rho ) = \vol(X,g_0)$.
\end{proposition}

\begin{preuve} On se propose d'appliquer le th{\'e}or{\`e}me 6.16 de \cite{Mar},
p.~332. On note $G = \ixg$, c'est un groupe alg{\'e}brique d{\'e}fini sur
$\rrf$ et semi-simple. Pour utiliser le r\'esultat 6.16 de \cite{Mar} il faut travailler avec des groupes de Lie
connexe, or $\Gamma $ est un sous-groupe de $G_+$, le sous-groupe de
$G$ constitu{\'e} des isom{\'e}tries pr{\'e}servant l'orientation et $G_+$ n'est
pas n{\'e}cessairement connexe. En effet, si $\gamma  = (\gamma _1\ld
\gamma _n)$, o{\`u} $\gamma _i \in \Isom(\wt X_i,\tilde g^i_0)$ , et si un
nombre pair de $\gamma _i$ renverse l'orientation alors $\gamma  \in
G_+$, n{\'e}anmoins $\gamma $ ne peut pas {\^e}tre connect{\'e} {\`a} l'identit{\'e}.

On rappelle que $G$ {\`a} un nombre fini de composantes connexes car c'est
un groupe alg{\'e}brique. Soit  $G^0$ la composante de l'{\'e}l{\'e}ment neutre et
$\Gamma ^0 = G^0 \cap \Gamma $.

Il est ais{\'e} de v{\'e}rifier que $G^0 = \Isom_+(\wt X_1,\tilde g^1_0)\times
\cdots \times \Isom_+(\wt X_p,\tilde g^p_0)$ o{\`u} $\Isom_+$ d{\'e}signe le
groupe (connexe) d'isom{\'e}tries directes.

Les quatre lemmes qui suivent n'utilisent pas l'irr\'eductibilit\'e de $\Gamma$. Cette hypoth\`ese ne sera utilis\'ee que pour appliquer le th\'eor\`eme de super-rigidit\'e.
\begin{lemme} 
Le groupe $\Gamma ^0$ est un r{\'e}seau cocompact de $G^0$
ainsi que de $G$.
\end{lemme}

\begin{preuve} L'application naturelle $\Gamma /\Gamma ^0
\hookrightarrow G/G^0$ est injective, $ \Gamma^0$ est donc d'indice fini
dans $\Gamma $ et est un r{\'e}seau cocompact de
$G$. Par ailleurs, $G^0/\Gamma ^0$  est une composante connexe de
$G/\Gamma ^0$, donc est compacte. Un th{\'e}or{\`e}me g{\'e}n{\'e}ral est prouv{\'e} dans
\cite{Rag}, p.~23 (th{\'e}or{\`e}me~1.13).
\end{preuve}

Pour all\'eger les notations nous d\'esignerons maintenant par $\rho$ la repr{\'e}sentation restreinte {\`a} $\Gamma
^0$. Soit $\Gamma ^1 = \rho ^{-1} \bigl(\rho (\Gamma ^0)\cap G^0\bigr)$.

\begin{lemme}
Le groupe $\Gamma ^1$ est d'indice fini dans $\Gamma^0$.
\end{lemme}

\begin{preuve} L'application $\Gamma ^0/\Gamma ^1 \to
G/G^0$ induite par $\rho $ est injective, d'o{\`u} le
r{\'e}sultat.
\end{preuve}

Le groupe $\Gamma ^1$ est donc un r{\'e}seau cocompact de $G$ (et de
$G^0$) qui de plus, comme $\Gamma $, est irr{\'e}ductible. La restriction
de $\rho $ {\`a} $\Gamma ^1$ est un homomorphisme
$$\rho  : \Gamma ^1 \to G^0$$
{\`a} valeurs dans le groupe semi-simple, connexe $G^0$.

\begin{lemme}
Les groupes $G$ et $G^0$ n'ont pas de centre.
\end{lemme}

\begin{preuve} Si $a\in G^0$ est dans le centre de $G^0$,
$a$ doit commuter avec tous les {\'e}l{\'e}ments de $G^0$~; or, pour $x \in
\wt X$ fix{\'e} $x = (x_1\ld x_n)$ les isom{\'e}tries du type $(\gamma _1\ld
\gamma _n)$, o{\`u} $\gamma _i$ est une isom{\'e}trie directe fixant $x_i$,
sont dans $G^0$. L'{\'e}l{\'e}ment $a$ doit donc fixer $x$, pour tout $x$,
c'est donc l'identit{\'e}.
\end{preuve}

Dans la terminologie de \cite{Mar}, le groupe $G^0$ est adjoint (il n'a pas
de centre et est d{\'e}fini sur $\rrf$, voir \cite{Mar}, p.~13).

\begin{lemme}\label{zariskidense}
Le groupe $\rho (\Gamma ^1)$ est Zariski-dense dans
$G^0$.
\end{lemme}

\begin{preuve} Rappelons que $\Gamma ^1$ est d'indice fini
dans $\Gamma $, on voit alors, de mani{\`e}re {\'e}l{\'e}mentaire, que
$$\vol(\rho ) = [\Gamma :\Gamma ^1]\vol (\rho_{| \Gamma ^1})$$
de sorte que l'hypoth{\`e}se de la proposition \ref{4.1} implique que $\vol(\rho
_{| \Gamma ^1})\ne 0$. La proposition \ref{reductif}, qui est un corollaire
de la remarque~1.4~{\it i)\/} de \cite{Lab}, montre que $\rho (\Gamma ^1)$
est r{\'e}ductif. Soit $H$ son adh{\'e}rence de Zariski, alors $H$ est
{\'e}galement r{\'e}ductif. Comme $H$ est alg{\'e}brique, quitte {\`a} restreindre {\`a}
un r{\'e}seau d'indice fini dans $\Gamma ^1$, on peut supposer que $H$ est
connexe.

L'alg{\`e}bre de Lie de $H$, \cad $\hg$ est une sous-alg{\`e}bre r{\'e}ductive
alg{\'e}brique de $\ggg^0$, alors d'apr{\`e}s le th{\'e}or{\`e}me~4, p.~261 de \cite{Oni-Ver}, il
existe une involution de Cartan de $\ggg^0$ qui stabilise $\hg$. Plus
pr{\'e}cis{\'e}ment, l'espace sym{\'e}trique $\wt X$ est identifi{\'e} {\`a} $G^0/K^0$ par
le choix d'une d{\'e}composition de Cartan de $G^0$ (ici, $K^0$ d{\'e}signe un
sous-groupe compact maximal de $G^0$)~; alors, si $\sigma ^0$ d{\'e}signe
l'involution de Cartan correspondante, il existe $g \in G^0$ tel que
l'involution  $g\sigma ^0 g^{-1}$ pr{\'e}serve $h$. Soit $x \in \wt X$
le point correspondant {\`a} la classe de $g$ dans $G^0/K^0$, alors
d'apr{\`e}s la proposition~2.6.2 de \cite{Ebe}, la sous-vari{\'e}t{\'e} $Hx = Y$ est
totalement g{\'e}od{\'e}sique dans $(\wt X,g_0)$ et $\rho (\Gamma ^1)$
invariante. 

On peut donc, pour calculer le volume de la repr{\'e}sentation $\rho
_{|\Gamma ^1}$ ($\rho $ restreinte {\`a} $\Gamma ^1$) utiliser une
application $C^1$, $\rho _{| \Gamma ^1}$-{\'e}quivariante de $\wt X$
dans $Y$. Si $H \ne G^0$ alors $\dim Y < \dim \wt X$ ce qui implique $\vol(\rho
_{| \Gamma ^1})=0$ et $\vol (\rho ) =0$.
Ceci est en contradiction avec l'hypoth{\`e}se de la
proposition.
\end{preuve}

Nous sommes maintenant en situation pour appliquer le th{\'e}or{\`e}me de
super-rigidit{\'e} 6.16~{\it b)\/} de \cite{Mar}, p.~332 (le groupe $G^0$, qui
est le groupe de d{\'e}part et d'arriv{\'e}e n'a aucune composante simple
compacte, \cad n'a pas de facteur $\rrf$-anisotrope). La repr{\'e}sentation
$\rho $ se prolonge en un (unique) homomorphisme continu
$$\tilde \rho  : G^0 \to G^0$$
qui est donc analytique (\cite{Hel}, p.~117, th{\'e}or{\`e}me~2.6). Le noyau $\Ker
\tilde \rho $ est un sous-groupe de Lie de $G^0$ (car ferm{\'e}). On
rappelle que $G^0 = \prodd^p_{i=1} G_i$, o{\`u} $G_i = \Isom_+(\wt
X_i,\tilde g^i_0)$ est un groupe simple. Comme $\Ker \tilde \rho $ est
normal, il est produit de certains $G_i$ de la liste pr{\'e}c{\'e}dente : $\Ker\tilde \rho  = \prodd^q_{k=1} G_k$ 
pour $q \le p$.

De plus, l'image $\tilde \rho (G^0)$ est un groupe de Lie isomorphe {\`a}
$G^0/\Ker \tilde \rho $, \cad isomorphe {\`a}  $\prodd^p_{k=q+1} G_k$ si
$q < p$, {\`a} $\{e\}$ sinon. En particulier $\tilde \rho (G^0)$ est un
sous-groupe semi-simple de $G^0$, invariant par $\rho (\Gamma
^1)$. L'argument du lemme \ref{zariskidense} montre que l'hypoth{\`e}se $\vol(\rho ) \ne
0$ implique que $\tilde \rho (G^0)$ doit {\^e}tre {\'e}gal {\`a} $G^0$, \cad que
$\tilde \rho $ doit {\^e}tre un automorphisme (analytique). En particulier
$\tilde \rho $ est un diff{\'e}omorphisme et $\tilde \rho (\Gamma ^1)$ est
un groupe discret et cocompact. Le th{\'e}or{\`e}me de Mostow permet de
conclure que les vari{\'e}t{\'e}s localement sym{\'e}triques $\wt X/\Gamma ^1$ et
$\wt X/\tilde \rho (\Gamma ^1)$  sont isom{\'e}triques et donc que
$$\vol (\rho _{|\Gamma ^1})= \vol (\wt X/\Gamma ^1)~.$$

Comme $\Gamma ^1$ est d'indice fini dans $\Gamma $, on en d{\'e}duit que
$$\vol(\rho ) = \vol(X,g_0)~.$$
\end{preuve}

Nous allons maintenant \'etudier les cas o\`u le r\'eseau $\Gamma$ est r\'eductible. 
On rappelle qu'un r{\'e}seau $\Gamma ^1$
de $G^0$ qui est r\'eductible v\'erifie les propri\'et\'es suivantes (voir \cite{Rag}, p.~86, 5.22) : il existe une famille finie de sous-groupes normaux
et connexes de $G^0$, $H_1\ld H_k$ telle que~:
\begin{itemize}
\item[i)] $H_i \cap \prodd_{j\ne i} H_j$ est discret pour
tout $i \in \{1\ld k\}$.

\item[ii)] $G^0 = \prodd^k_{i=1} H_i$.

\item[iii)] $\Gamma ^1_i = H_i \cap \Gamma ^1$ est un r{\'e}seau
irr{\'e}ductible de $H_i$.

\item[iv)] $\prodd^k_{i=1} \Gamma ^1_i$ est un sous-groupe
normal d'indice fini de $\Gamma ^1$.
\end{itemize}

Comme pr{\'e}c{\'e}demment nous pouvons travailler {\`a} un sous-groupe d'indice
fini pr{\`e}s et donc supposer que $\prodd^k_{i=1} \Gamma ^1_i = \Gamma
^1$. De m{\^e}me, chaque $H_i$ doit {\^e}tre un produit de facteurs simples
composant $G^0$, \cad
$$H_i = \prodd^{p_i+r_i}_{s=p_i} G_s~.$$

En particulier, si  $i\ne j$, $H_i$ et $H_j$ commutent.

\begin{proposition}\label{reductible}
Avec les notations ci-dessus, si $\Gamma$ est r\'eductible et $\vol(\rho ) \ne
0$ et si, pour tout $i$, $H_i$ est super-rigide alors $\rho (\Gamma )$
est un r{\'e}seau cocompact et donc $\vol(\rho ) = \vol(X,g_0)$.
\end{proposition}

\begin{preuve} Par groupe super-rigide nous entendons un groupe auquel nous
pouvons appliquer le th{\'e}or{\`e}me de super-rigidit{\'e}, \cad, dans notre
situation, soit $H_i$ est de rang sup{\'e}rieur o{\`u} {\'e}gal {\`a} 2 ($r_i\ge 1$)
ou bien $H_i$ est le groupe d'isom{\'e}tries directes d'un espace
hyperbolique quaternionien ou du plan hyperbolique de Cayley.

Comme pr{\'e}c{\'e}demment nous travaillons avec le sous-groupe $\Gamma
^1$. D{\'e}finissons $\rho _i = \rho _{|\Gamma ^1_i}$, pour $i=1\ld k$~;
ici nous commettons un abus de langage et identifions $\Gamma ^1_i$ et
$\{e\} \times \cdots \times \{e\} \times \Gamma ^1_i \times \{e\}
\times \cdots \times \{e\}$. D{\'e}finissons les groupes $K_i =
\overline{\rho_i (\Gamma ^1_i)}$, l'adh{\'e}rence de Zariski de $\rho
_i(\Gamma ^1_i)$~; ce sont des sous-groupes alg{\'e}briques de $G^0$ et le
th{\'e}or{\`e}me 6.15~{\it i), a)\/} de \cite{Mar}, p.~332 affirme que, si rang
$H_i\ge 2$, $K_i$ est un groupe semi-simple. Insistons sur le fait que $H_i$ est consid\'er\'e comme un sous-groupe du groupe de d\'epart de la repr\'esentation $\rho$ et $K_i$ comme un sous-groupe du groupe d'arriv\'ee. Le m{\^e}me r{\'e}sultat pour le cas o{\`u} $H_i$ est le groupe d'isom{\'e}tries (directes) de l'espace
hyperbolique quaternionien ou du plan hyperbolique de Cayley est
prouv{\'e} dans \cite{Sta}. Dans tous les cas, donc, $K_i$ est un groupe
semi-simple.

\begin{lemme}\label{5.8}
Les groupes $K_i$ sont sans facteurs compacts et $\prod K_i$ est semi-simple.
\end{lemme}

\begin{preuve} Quitte {\`a} passer {\`a} un sous-groupe d'indice fini de $\Gamma
^1_i$ nous pouvons supposer que $K_i$ est un produit de groupes
simples. Supposons qu'il contienne un facteur compact, soit
$$K_i = L^1_i \times \cdots \times L^{k_i}_i \times U$$
o{\`u} $U$ est un groupe simple compact~; alors $U$ est normal dans
$K_i$. Par ailleurs $\rho _i(\Gamma ^1_i)$ et $\rho _j(\Gamma ^1_j)$
commutent si $i\ne j$, car $\Gamma ^1_i$ et $\Gamma ^1_j$ commutent
dans $G^0$ (nous faisons ici l'abus de langage signal{\'e} pr{\'e}c{\'e}demment)~;
les groupes $K_i$ et $K_j$ commutent donc {\'e}galement si $j\ne i$. Le
groupe $U$ est donc normal dans le produit $\prodd^k_{i=1} K_i$; remarquons que ce
produit est d{\'e}fini comme le groupe engendr\'e par les produits d'\'el\'ements de $K_i$. Il est d\'efini sans ambigu{\"\i}t{\'e} car les groupes commutent deux {\`a} deux. Par ailleurs $\prodd^k_{i=1} K_i$ est Zariski-dense donc \'egal \`a $G_0$; en effet, il contient $\rho (\Gamma ^1)$ qui Zariski dense car $\vol(\rho ) \ne 0$ (voir le lemme \ref{zariskidense}). Le groupe $U$ est donc normal dans $G^0$~; il est alors {\'e}gal {\`a} une des composantes de $G^0$ ou bien r{\'e}duit {\`a} l'{\'e}l{\'e}ment neutre~; aucune des composantes de $G^0$ n'{\'e}tant compacte $U$ est
trivial.

Enfin les arguments pr\'ec\'edents montrent que les composantes simples de $K_i$ et $K_j$ sont distincts et donc que $\prodd^k_{i=1} K_i$ est semi-simple (c'est le produit de toutes les composantes des groupes $K_i$).
\end{preuve}

Nous pouvons donc appliquer le th{\'e}or{\`e}me de super-rigidit{\'e} de \cite{Mar},
(6.16 {\it c)}, p.~332) pour les composantes $H_i$ de rang $\ge 2$ et
celui de \cite{Cor} pour les autres et affirmer que les repr{\'e}sentations
$\rho _i$
se prolongent en des morphismes continus
$$\varphi _i : H_i \to K_i~.$$
On construit alors un prolongement de $\rho $ en
\begin{align*}
\varphi  : G^0 &\la \prodd^k_{i=1}K_i\subset G^0\\
(\gamma _1\ld \gamma _k )  &\longmapsto \big(\varphi _1(\gamma
_1),\cdots ,\varphi _k (\gamma _k)\big)
\end{align*}
Les morphismes $\varphi _i$ commutent et $\varphi $ est bien d{\'e}fini et
est un morphisme continu. On termine donc la preuve de la proposition ref{reductible} par les m{\^e}mes
arguments que ceux de la preuve de la proposition \ref{irreductible}.
\end{preuve}

Nous nous int\'eressons maintenant au cas o\`u $G_0$ poss\`ede des composantes simples non super-rigides.
Supposons donc que $\Gamma ^1 = \Gamma ^1_1 * \Gamma ^1_2$ (produit libre) o{\`u}
$\Gamma ^1_1$ est un r{\'e}seau cocompact d'un groupe $H_1$ extension finie d'un produit de
groupes super-rigides et $\Gamma ^1_2$ est un
r{\'e}seau cocompact de $H_2$ produit de copies de $PO(k,1)$ et
$PU(k',1)$. Les arguments qui pr{\'e}c{\`e}dent s'appliquent pour montrer que

{\it i)\/} $\rho (\Gamma )$ est Zariski dense dans $G^0$ si $\vol(\rho
) \ne 0$.

{\it ii)\/} Soit $K_i$, $i=1,2$, l'adh{\'e}rence de Zariski de $\rho
_i(\Gamma ^1_i)$. La densit{\'e} de $\rho (\Gamma )$ implique la densit{\'e}
(pour la topologie de Zariski) de $K_1K_2$~; les deux groupes
commutent. En d{\'e}composant $\Gamma ^1_1$ en produits de r{\'e}seaux
cocompacts irr{\'e}ductibles on voit, en utilisant les arguments de la
preuve du lemme \ref{5.8}, que $K_1$ est semi-simple sans facteurs
compacts. Si $G^0 = \prodd^p_{i=1} G_k$, o{\`u} les $G_k$ sont des groupes
d'isom{\'e}tries directes d'espaces sym{\'e}triques de rang 1 et de type non
compact, alors $K_1 = \prodd^q_{k =1} G_k$ (par exemple)~; en effet,
$K_1$ est un sous-groupe normal de $G^0$. Le  groupe $K_2$ est donc
inclus dans $\prodd^p_{q+1} G_k$ et comme $K_1K_2$ est Zariski dense,
$K_2 = \prodd^p_{q+1} G_k$ (en particulier il est semi-simple). On
peut choisir une application {\'e}quivariante $f_i$ de $H_i \to K_i$ et un
calcul imm{\'e}diat montre que
$$\vol(\rho )= \vol (\rho _1) \vol(\rho _2)~,$$
de sorte que $\vol(\rho _1) \ne 0$. Le th{\'e}or{\`e}me de super-rigidit{\'e}
(appliqu{\'e} comme pr{\'e}c{\'e}demment aux composantes irr{\'e}ductibles de $\Gamma
^1_1$) permet d'{\'e}tendre $\rho _1$ en un morphisme continu
$$\varphi _1 : H_1 \to K_1$$
et la non nullit{\'e} de $\vol(\rho _1)$ montre que $\varphi _1$ est un
isomorphisme et donc montre que $\rho (\Gamma ^1_1)$ est isomorphe {\`a}
$\Gamma ^1_1$ ce qui conduit {\`a}
$$\vol (\rho _1) = \vol (H_1/\Gamma ^1_1)~.$$
Le groupe $K_2$ contient les composantes non super-rigides de $G^0$,
de m{\^e}me que $H_2$, ils sont donc isomorphes. C'est la seule composante
non triviale de $\rho $.

Nous donnons maintenant un exemple de repr{\'e}sentation du groupe
fondamental d'une vari{\'e}t{\'e} hyperbolique r{\'e}elle, de volume non nul et
d'image non discr{\`e}te.

\begin{exemple} : \textbf{produit amalgam{\'e}.}
\end{exemple}
Soit $X$ une vari{\'e}t{\'e} hyperbolique compacte de dimension $n\ge
3$. Supposons qu'il existe dans $X$ une hypersurface compacte plong{\'e}e
totalement g{\'e}od{\'e}sique not{\'e}e $\Sigma $ et incompressible, \cad telle
que l'application induite~: $\pi _1(\Sigma ) \to \Pi _1(X)$ soit une
injection. Nous supposons de plus que cette hypersurface s{\'e}pare $X$ en
deux composantes connexes $X_A$ et $X_B$ de groupe fondamental
respectif $A$ et $B$. En posons $C = \pi _1(\Sigma )$, le th{\'e}or{\`e}me de
Van Kampen montre que
$$\pi _1(X) = A *_C B$$
produit amalgam{\'e} de $A$ et $B$ sur $C$. Les groupes $\pi _1(X)$, $A$,
$B$ et $C$ sont des sous-groupes de $PO(n,1)$ et agissent donc sur
l'espace hyperbolique $\hhf^n$. Choisissons un relev{\'e} $\wt \Sigma $ de
$\Sigma $ dans $\hhf^n$~; $\wt \Sigma $ est une hypersurface totalement
g{\'e}od{\'e}sique. On identifie $C$ au sous-groupe de $\pi _1(X)$ qui fixe
$\wt \Sigma $. Soit $s$ la sym{\'e}trie par rapport {\`a} $\wt \Sigma $,  on
d{\'e}finit
\begin{align*}
\rho  : \pi _1(X) &\la PO(n,1)\\
a\in A &\longmapsto a\\
b\in B &\longmapsto sbs^{-1}~.\end{align*}

\begin{lemme}
L'application $\rho $ d{\'e}finit une repr{\'e}sentation de
$\pi _1(X)$ dans $PO(n,1)$.
\end{lemme}

\begin{preuve} Le groupe $\pi _1(X)$ est le quotient du produit libre $A*B$
par les relations qui consistent {\`a} identifier un {\'e}l{\'e}ment de $C$ dans
$A$ avec le m{\^e}me {\'e}l{\'e}ment dans $B$. Comme $ \rho $ est un morphisme en
restriction {\`a} $A$ et {\`a} $B$ respectivement, il suffit de v{\'e}rifier la
compatibilit{\'e} avec les relations. Or, si $c\in C$ 
$$s c s^{-1} = c$$
d'o{\`u} le r{\'e}sultat.
\end{preuve}

Afin de calculer le volume de cette repr{\'e}sentation il faut trouver une
application lipschitzienne $f : \hhf^n \to \hhf^m$, $\rho
$-{\'e}quivariante.

\begin{proposition}
Avec les notations ci-dessus on a,
$$\vol(\rho ) = \vol(X_A) - \vol (X_B)~.$$
\end{proposition}

\begin{preuve} Nous allons d{\'e}crire $f$ de mani{\`e}re pr{\'e}cise et le calcul du
volume s'ensuivra. Le fait que $\pi _1(X)$ soit un produit amalgam{\'e}
est {\'e}quivalent (\cite{Ser}, p.~48) {\`a} l'existence d'un arbre $T$, sur lequel
$\pi _1(X)$ op{\`e}re (sans inversion) en sorte que le quotient soit un
segment (deux sommets joints par une ar{\^e}te). Les sous-groupes $A$,
$B$ et $C$ sont alors les stabilisateurs respectifs des deux sommets
et de l'ar{\^e}te de l'arbre quotient. Nous allons donner une description
g{\'e}om{\'e}trique de cet arbre $T$. Nous avons choisi un relev{\'e} $\wt \Sigma
$ de l'hypersurface compacte $\Sigma $ plong{\'e}e dans $X$~; $\Sigma $
{\'e}tant une sous-vari{\'e}t{\'e} plong{\'e}e, sans auto-intersection, les translat{\'e}s
$\gamma  \wt \Sigma $ de $\wt \Sigma $ par les {\'e}l{\'e}ments $\gamma  \in
\pi _1(X)$ sont deux {\`a} deux disjoints~; ils s{\'e}parent donc $\hhf^n$ en
une infinit{\'e} de composantes connexes. Les deux composantes connexes
dont l'adh{\'e}rence contient $\wt \Sigma $ sont des rev{\^e}tements
universels de $X_A$ et $X_B$ respectivement, que nous noterons $\wt
X_A$ et $\wt X_B$. Les autres composantes connexes sont les translat{\'e}s
par les {\'e}l{\'e}ments de $\pi _1(X)$ de $\wt X_A$ et $\wt X_B$. Les
sous-groupes $A$ et $B$ pr{\'e}servent $\wt X_A$ et $\wt X_B$
respectivement (apr{\`e}s un choix convenable d'un point base et d'un de
ses relev{\'e}s).

Maintenant, choisissons un point $x_a \in \wt X_A$ et un point $x_b
\in \wt X_B$, les sommets de l'arbre $T$ sont les $\gamma (x_a)$ et
$\gamma  (x_b)$, o{\`u} $\gamma $ parcourt $\pi _1(X)$~; on joint deux
sommets $\gamma  (x_a)$ et $\gamma '(x_b)$ (de type diff{\'e}rent) si, et
seulement si, les composantes connexes correspondantes $\gamma \wt
X_A$ et $\gamma ' \wt X_B$ sont telles que $\ol{\gamma \wt X_A} \cap
\ol {\gamma '\wt X_B} \ne \emptyset$.

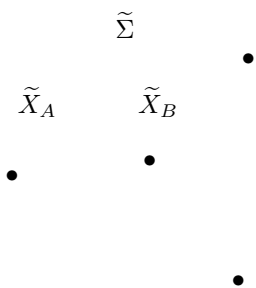
\InsertFig 40 0 {besson1.eps}
\LabelTeX 0 60 \ELTX
\LabelTeX 28 45 $\wt \Sigma $\ELTX
\LabelTeX 15 35 $\wt X_A$\ELTX
\LabelTeX 31 35 $\wt X_B$\ELTX
\LabelTeX 31.7 27.8 $\bullet $\ELTX 
\LabelTeX 13.4  25.8 $\bullet $\ELTX 
\LabelTeX 44.8 41.4 $\bullet $\ELTX 
\LabelTeX 43.5 11.8 $\bullet $\ELTX 
\EndFig 

Soit alors $x \in \wt X$, il appartient {\`a} une composante connexe du
compl{\'e}mentaire de $\cupp_{\gamma \in \Gamma } \gamma \wt \Sigma $ qui
correspond {\`a} un sommet  de l'arbre pr{\'e}c{\'e}dent. Dans cet arbre il existe
un unique chemin joignant la composante $\wt X_A$ {\`a} celle de $x$~; ce
chemin est une succession d'ar{\^e}tes $e_1,e_2\ld e_k$ prises dans
l'ordre, de la composante $\wt X_A$ {\`a} celle de $x$. Chacune de ces
ar{\^e}tes correspond {\`a} une image de $\wt \Sigma$ et nous noterons
$s_{e_i}$ la sym{\'e}trie orthogonale hyperbolique par rapport {\`a} cette
hypersurface totalement g{\'e}od{\'e}sique.

\begin{definition}
On pose $f(x) = s_{e_1} \circ s_{e_2} \circ
\cdots \circ s_{e_k}(x)$
\end{definition}

L'application $f$ est bien d{\'e}finie. Elle est $\ci$ par morceaux et
continue~; en effet, la seule ambigu{\"\i}t{\'e} dans la formule ci-dessus est
lorsque $x$ est sur l'hypersurface d{\'e}finie par $e_k$, mais dans ce cas
$s_{e_k}(x) = x$.

\begin{lemme}
L'application $f$ est $\rho $-{\'e}quivariante.
\end{lemme}

\begin{preuve} Il suffit de v{\'e}rifier l'{\'e}quivariance pour les {\'e}l{\'e}ments de $A$
et ceux de $B$ qui engendrent le groupe fondamental de $X$.

{\it a)\/} Si $a \in A$, le chemin dans l'arbre joignant la composante
$\wt X_A$ {\`a} celle de $ax$ est constitu{\'e} des ar{\^e}tes $ae_1 \ld
ae_k$~; en effet, puisque $a \in A$, $a\wt X_A = \wt X_A$ et
l'ar{\^e}te $ae_1$ a son origine dans $\wt X_A$. D'o{\`u}
\begin{align*}
f(ax) &= s_{ae_1} \cdots s_{ae_k}(ax)\\
&= as_{e_1} a^{-1} \cdots a s_{e_k}a^{-1}(ax)\\
&= a\ f(x) = \rho (a) f(x)\end{align*}

{\it b)\/} si $b\in B$, le chemin joignant la composante $\wt X_A$ {\`a}
celle de $bx$ est constitu{\'e} du chemin dans l'arbre joignant $\wt X_A$
{\`a} $b\wt X_A$ suivi de l'image par $b$ du chemin pr{\'e}c{\'e}dent. Rappelons
que les sommets de l'arbre sont les {\'e}l{\'e}ments de $\Gamma /A$ et $\Gamma
/B$ et les ar{\^e}tes sont les {\'e}l{\'e}ments de $\Gamma /C$ (voir \cite{Ser}). Par
exemple, la composante connexe $\wt X_A$ correspond {\`a} $eA$ (classe de
l'{\'e}l{\'e}ment neutre $e$), celle de $\wt X_B$ {\`a} $eB$~; elles sont

\InsertFig 40 0 {besson2.eps}
\LabelTeX -10 10 \ELTX
\LabelTeX 0 3 $eA $\ELTX
\LabelTeX 20  3 $eB$\ELTX
\LabelTeX 42 3 $bA$\ELTX
\LabelTeX 10 -3 $eC$\ELTX
\LabelTeX 30 -3 $bC$\ELTX
\LabelTeX 0.2 -0.2 $\bullet $\ELTX 
\LabelTeX 20.2  -0.2 $\bullet $\ELTX 
\LabelTeX 42.2 -0.2 $\bullet $\ELTX 
\EndFig

\noindent reli{\'e}es par l'ar{\^e}te $eC$.  Par ailleurs $b\wt X_A$
correspond {\`a} la classe $bA$ reli{\'e}e {\`a} $eB$ par l'ar{\^e}te $bC$. En
conclusion, nous avons
$$f(bx) = s \circ s_{be} \circ (bs_{e_1} \cdots s_{e_k} b^{-1})
(bx)$$
o{\`u} $e$ d{\'e}signe par abus de langage l'ar{\^e}te $eC$ et $s_e = s$. D'o{\`u}
\begin{align*}
f(bx) &= (sbsb^{-1}) \circ (bs_{e_1}\cdots s_{e_k} b^{-1}) (bx)\\
&= \rho (b) f(x)~.
\end{align*}
\end{preuve}

\vskip5pt
{\it Fin de la preuve de la proposition \pointir} La fin de la preuve est
{\'e}vidente~; en effet $f$ renverse l'orientation sur $\wt X_B$ et est
l'identit{\'e} sur $\wt X_A$, il suffit donc de choisir un domaine
fondamental dans la r{\'e}union $\wt X_A\cup \wt X_B$ pour lequel $\Sigma
$ se rel{\`e}ve sur $\wt \Sigma $.
\end{preuve}

Pour {\^e}tre complet, il faut construire des vari{\'e}t{\'e}s $X$ hyperboliques
admettant une hypersurface connexe s{\'e}parante qui s{\'e}pare la vari{\'e}t{\'e} en
deux parties de volume distinct. Cette construction nous a {\'e}t{\'e}
sugg{\'e}r{\'e}e par N.~Bergeron. Soit $M_1$ une vari{\'e}t{\'e} compacte de dimension
3, hyperbolique {\`a} bord totalement g{\'e}od{\'e}sique qui est une surface
compacte connexe not{\'e}e $\Sigma $. De tels exemples existent (voir
\cite{Ota}, \cite{Bon}, th.~4.3 et \cite{Ber}). Consid{\'e}rons le double $M$ obtenu par
recollement de deux copies de $M_1$ le long de $\Sigma $. La vari{\'e}t{\'e}
compacte $M$ sans bord est hyperbolique car $\Sigma $ est totalement
g{\'e}od{\'e}sique. Le th{\'e}or{\`e}me 2 de \cite{Ber2} montre que l'on peut construire un
rev{\^e}tement fini $\wh M$ de $M$ tel que $\Sigma $ se rel{\`e}ve
isom{\'e}triquement {\`a} $\wh M$ en une sous-vari{\'e}t{\'e} totalement g{\'e}od{\'e}sique
$\wh \Sigma $ non s{\'e}parante. On d{\'e}coupe alors $\wh M$ le long de $\wh
\Sigma $ pour obtenir une vari{\'e}t{\'e} {\`a} bord dont les deux composantes du
bord, not{\'e}e $\wh \Sigma _1$, $\wh \Sigma _2$, sont isom{\'e}triques {\`a} $\wh
\Sigma  \simeq \Sigma $ et on recolle {\`a} chacune de ces composantes une
copie de $M_1$. Alors, $\wh \Sigma _1$ (et $\wh \Sigma _2$) d{\'e}coupe la
nouvelle vari{\'e}t{\'e} hyperbolique en deux composantes l'une de volume {\'e}gal
{\`a} $\vol (M_1)$ et l'autre de volume {\'e}gal {\`a} $\vol(\wh M) + \vol(M_1) >
\vol(M_1)$.
\end{preuve}

Il serait int{\'e}ressant de disposer de tels exemples en dimension $n\ge
4$. Remarquons, par ailleurs, que l'ensemble des valeurs de $\vol
(\rho )$ ainsi obtenu est discret (pour une vari{\'e}t{\'e} donn{\'e}e)~; une
explication pr{\'e}cise {\`a} ce ph{\'e}nom{\`e}ne est fournie par le chapitre suivant.

\section{Volume et d{\'e}formations}\label{applications}

Nous avons d{\'e}j{\`a} remarqu{\'e} que, lorsque la dimension de $X$ est paire, le volume
d'une repr{\'e}sentation ($\widetilde X$ est suppos{\'e}e sym{\'e}trique) est le nombre d'Euler
du fibr{\'e} plat correspondant. En particulier, ce nombre est constant le long des 
d{\'e}formations continues de repr{\'e}sentations. Nous allons prouver un r{\'e}sultat analogue
dans le cas o{\`u} la dimension de $X$ est impaire. De telles d{\'e}formations existent en dimension
$3$ (\cite{Ben-Pet}) et nous en donnons des exemples. La constance du volume est prouv{\'e}e en
dimension $3$ par S.~Reznikov \cite{Rez2}; nous donnons ici une preuve, valable en
toute dimension, qui repose sur le formule de Schl{\"a}fli. Dans ce qui suit
$M$ d{\'e}signe une vari{\'e}t{\'e} riemannienne ferm{\'e}e et
orient{\'e}e de dimension $n$ et $\widetilde X$
l'espace hyperbolique r{\'e}el simplement connexe de dimension $n$.

\begin{theoreme}
\label{vol-const}
Soit $M$ une vari{\'e}t{\'e} diff{\'e}rentielle ferm{\'e}e et orient{\'e}e et 
$\rho_t : \Pi_1(M)\longrightarrow Isom({\widetilde X})$ une famille de repr{\'e}sentations
qui d{\'e}pend de mani{\`e}re $C^1$ du param{\`e}tre $t\in \bf{R}$, alors le volume 
$\vol (\rho _t)$ est constant.
\end{theoreme}

La preuve repose sur un lemme technique dont le but est de construire une application
{\'e}quivariante affine par morceaux particuli{\`e}re. Par application affine nous entendons une application affine le long de toute g\'eod\'esique.

\begin{lemme}\label{triangulation}
 Sous les hypoth{\`e}ses du th{\'e}or{\`e}me \ref{vol-const}, il existe une triangulation ${\cal T}$
de $\widetilde M$ et une application continue et affine par morceaux $\tilde {f_0} : {\wt M}\to {\wt X}$ qui
est $\rho _0$-{\'e}quivariante et non d{\'e}g{\'e}n{\'e}r{\'e}e au sens o{\`u} l'image par $\tilde {f_0}$ d'un simplexe
de la triangulation ${\cal T}$ est un simplexe g{\'e}od{\'e}sique de $\widetilde X$ non d{\'e}g{\'e}n{\'e}r{\'e}.
\end{lemme}

\begin{preuve}
Un th{\'e}or{\`e}me classique affirme que toute vari{\'e}t{\'e} lisse compacte $M$ est hom{\'e}omorphe
{\`a} un complexe simplicial $K$; plus pr{\'e}cis{\'e}ment $K$ est un espace triangul{\'e} muni d'une m\'etrique euclidienne par morceaux (que l'on peut r{\'e}aliser dans $\bf R^n$). Cet hom{\'e}omorphisme peut, de plus, {\^e}tre choisi Lipschitzien. Le volume de toute repr{\'e}sentation de $\Pi_1(M)=\Pi_1(K)$ peut donc se calculer en int{\'e}grant sur $M$
ou bien sur $K$. Dans la suite nous noterons {\'e}galement $M$ ce complexe euclidien par morceaux et
toute triangulation sera une subdivision de la d{\'e}composition de $K$ en simplexes.

Choisissons alors une triangulation suffisamment fine de $M$ et
appelons ${\cal T}_{\widetilde M}$
la triangulation invariante par $\Pi_1(M)$  sur $\widetilde M$ qui s'en d{\'e}duit par image r{\'e}ciproque.
Soit $D$ un domaine fondamental (de Dirichlet) dans $\widetilde M$ pour l'action de $\Pi_1(M)$.
Quitte {\`a} modifier un peu ${\cal T}_{\widetilde M}$ ou bien $D$ on peut supposer qu'aucun sommet de
la triangulation n'est sur $\partial D$.

Notons $(m_1,...,m_N)$ la liste des sommets de ${\cal T}_{\widetilde
M}$ qui sont dans l'int{\'e}rieur de $D$, N est alors le cardinal des
sommets de la triangulation de d{\'e}part sur $M$. Choisissons
maintenant $N$ points dans $\widetilde X$, not{\'e}s $(y_1,...,y_N)$
de sorte que si $(m_{i_1},...,m_{i_{k+1}})$ est un $k$-simplexe de
${\cal T}$ alors le simplexe g{\'e}od{\'e}sique de $\widetilde X$ de
sommets $(y_{i_1},...,y_{i_{k+1}})$ est non d{\'e}g{\'e}n{\'e}r{\'e}
pour tout $k\le \textrm{dim}M+1$.  Ceci est toujours possible car, pour chaque sommet
$y_j$, la r{\'e}union des conditions de d{\'e}g{\'e}n{\'e}rescence des
simplexes contenant $y_j$ est un ensemble ferm{\'e} d'int{\'e}rieur
vide (une r{\'e}union finie de $k-1$-plans). Ces choix {\'e}tant fait,
il existe autour de chaque point $y_j$ un petit voisinage $V_j$ en
sorte que, pour n'importe quel choix de points $y'_1,...,y'_N$ avec
$y'_j\in V_j$, la propri{\'e}t{\'e} de non d{\'e}g{\'e}n{\'e}rescence
ci-dessus soit encore v{\'e}rifi{\'e}e. Par la suite nous aurons \'egalement
besoin de choisir les point $y_j$ de sorte que
$$\forall \gamma\in \Pi_1(M),\, \forall j\ne i,\quad y_i\ne \rho_0
(\gamma)y_j\,,$$ 
ceci est toujours possible car la r{\'e}union des
points de l'orbite des $y_j$, pour $j\ne i$, qui sont dans $V_i$ est un
ensemble d{\'e}nombrable. On proc{\`e}de dons par r{\'e}currence,
$y_1$ {\'e}tant fix{\'e} on choisit $y_2\in V_2$ dans l'ensemble
partout dense qui est le compl{\'e}mentaire de l'orbite de $y_1$ ,
puis $y_3$ dans le compl{\'e}mentaire des orbites de $y_1$ et $y_2$ et
ainsi de suite.

On d{\'e}finit alors $\tilde {f_0}$ par :
$$\forall \gamma\in\Pi_(M),\,\forall i=1,...,N\quad \tilde{f_0}(\gamma
m_i)=\rho _0(\gamma )y_i\, ,$$ 
et on {\'e}tend $\tilde {f_0}$ {\`a}
l'int{\'e}rieur d'un simplexe $(m_{i_1},...,m_{i_N})$ en une
application affine sur le simplexe g{\'e}od{\'e}sique engendr{\'e} par
les points $y_{i_1},...,y_{i_N}$; on utilise pour cela la m{\'e}trique
euclidienne sur les simplexes de $\widetilde M$ et la m{\'e}trique
hyperbolique sur ceux de $\widetilde X$.  Par le choix des points
$y_i\in{\widetilde X}$, tous les simplexes dont les sommets sont dans
l'int{\'e}rieur de $D$ sont transform{\'e}s par $\tilde {f_0}$ en des
simplexes non d{\'e}g{\'e}n{\'e}r{\'e}s. Consid{\'e}rons 
 maintenant
le cas o{\`u} certains sommets sont dans l'int{\'e}rieur de $D$ et
d'autres {\`a} l'ext{\'e}rieur. Soit
$(m_{i_1},...,m_{i_p},\gamma_{j_1}m_{j_1},...,\gamma_{j_q}m_{j_q})$ un
tel simplexe et supposons que son image par $\tilde {f_0}$,
c'est-{\`a}-dire le simplexe not\'e
 $(y_{i_1},...,y_{i_p},\rho
(\gamma_{j_1})y_{j_1},..., \rho (\gamma_{j_q})y_{j_q})$, soit
d{\'e}g{\'e}n{\'e}r{\'e}; cela signifie qu'il existe $1\leq k\leq q$
tel que $\rho (\gamma_{j_k})y_{j_k}$ appartienne au sous-espace
totalement g{\'e}od{\'e}sique $E$ engendr{\'e} par les points
$(y_{i_1},...,y_{i_p},\rho (\gamma_{j_1})y_{j_1},...,\rho
(\gamma_{j_{k-1}})y_{j_{k-1}})$ (rappelons que, par construction,
$(y_{i_1},...,y_{i_p})$ est un $(p-1)$-simplexe non
d{\'e}g{\'e}n{\'e}r{\'e}); on d{\'e}place alors $y_{j_k}$ {\`a}
l'int{\'e}rieur de $V_{j_k}$ pour le s{\'e}parer de $\rho
(\gamma_{j_k})^{-1}E$; ceci est possible si $E$ reste fixe lorsque
l'on d{\'e}place $y_{j_k}$, c'est-{\`a}-dire si aucun des points
$y_{i_1},..., y_{i_p},\rho (\gamma_{j_1})y_{j_1},\\ ...,\rho
(\gamma_{j_{k-1}})y_{j_{k-1}}$ n'est dans l'orbite de $y_{j_k}$. Par
le choix des $y_i$ ceci ne peut se produire que si
$y_{j_k}=y_{j_l}\textrm{ avec } 
 l=1,...,q\textrm{ et }l\ne k$ ou
bien $y_{i_k}=y_{i_l},\, l=1,...,p$. Les points $y_i$ {\'e}tant en
bijection avec les points $m_i$ cela impliquerait que dans le simplexe
$(m_{i_1},...,m_{i_p},\gamma_{j_1}m_{j_1},...,\gamma_{j_q}m_{j_q})$
deux des points $m_l$ co{\"\i}ncident et donc qu'au quotient sur $M$
il se projette sur un simplexe d{\'e}g{\'e}n{\'e}r{\'e} ce qui est
impossible. On peut donc s{\'e}parer $y_{j_k}$ du sous-espace
totalement g{\'e}od{\'e}sique $\rho (\gamma_{j_k})^{-1}E$. On utilise
ensuite l'argument de densit{\'e} pour choisir le nouveau point
$y_{j_k}$ disjoint de la r{\'e}union des orbites par $\rho (\Pi_1(M))$
des autres points $y_l$. On proc{\`e}de alors par r{\'e}currence sur
les simplexes consid{\'e}r{\'e}s qui sont en nombre fini.

Les autres simplexes sont des images par un {\'e}l{\'e}ment $\rho
(\gamma )$, pour $\gamma\in\Pi_1(M)$ des simplexes d'un des deux types
pr{\'e}c{\'e}dents. Ceci prouve le lemme \ref{triangulation}.
\end{preuve}

 \textbf{preuve du th\'eor\`eme}

 Nous noterons
$\mathcal{T}_{\widetilde X}$ la collection des simplexes de
$\widetilde X$ ainsi obtenue. 
 Soit $F$ une face de codimension $2$
de $\mathcal{T}_{\widetilde M}$ et $F'$ son image dans 
$\mathcal{T}_{\widetilde X}$. L'\'etoile de $F$ dans
$\mathcal{T}_{\widetilde M}$ contient un nombre fini de $n$-simplexes
$s_1,...,s_k$ dont les images sont not\'ees $s'_1,...,s'_k$. Le link
autour de $F$ est un cercle. Pr\'ecis\'ement, consid\'erons un
voisinage tubulaire de rayon assez petit, not\'e Tub$(F)$, de cette
face $F$ de codimension $2$. Alors le bord de Tub$(F)$ est
diff\'eomorphe \`a $F\times S^1$. La vari\'et\'e $M$ est suppos\'ee
orient\'ee, et donc aussi $\widetilde M$. Sur le bord de Tub$(F)$ nous
choisissons une courbe $\mathcal{C}$ g\'en\'erateur de 
$H_1(\partial\textrm{Tub}(F),\textbf{Z})\simeq \textbf{Z}$; nous
pouvons, par exemple, prendre l'intersection de
$\partial\textrm{Tub}(F)$ avec un hyperplan orthogonal \`a $F$ en un
point (on peut d\'efinir un tel hyperplan bien que la m\'etrique sur
$\widetilde M$, qui est euclidienne sur chaque simplexe, soit
singuli\`ere en $F$). Si nous choisissons arbitrairement une
orientation sur chaque face de codimension $2$, donc en particulier
sur $F$, cela fournit une orientation du cercle $\mathcal{C}$
compatible avec celle de $\widetilde M$.
 
 L'application
$\tilde{f_0}$, lin\'eaire par morceaux, envoie $F$ sur $F'$ (par
construction) et donc $\partial\textrm{Tub}(F)$ sur un cylindre
topologique que l'on peut projeter, \`a partir de $F'$, sur le bord
$\partial\textrm{Tub}(F')$ d'un petit voisinage tubulaire de $F'$
(pour la m\'etrique hyperbolique). Cela induit une application,

$$\tilde{f_0}\,: H_1(\partial\textrm{Tub}(F),\textbf{Z})\longmapsto
H_1(\partial\textrm{Tub}(F'),\textbf{Z})$$

 et on appelle degr\'e transverse de $\tilde{f_0}$ en $F$, l'image
par $\tilde{f_0}_*$ du g\'en\'erateur de
$H_1(\partial\textrm{Tub}(F),\textbf{R})$; cette classe est un multiple entier de la classe fondamentale de $H_1(\partial\textrm{Tub}(F'),\textbf{Z})$ et nous pouvons donc, par abus de langage, identifi\'e le degr\'e transverse \`a un nombre entier relatif. On peut \'egalement
d\'efinir ce degr\'e en utilisant le cercle $\mathcal{C}$ trac\'e sur
$\partial\textrm{Tub}(F)$ et un cercle $\mathcal{C}'$ analogue sur
$\partial\textrm{Tub}(F')$ sur lequel on
 projette
$\tilde{f_0}(\mathcal{C})$.

Soit $\theta (F,s)$ (resp. $\theta '(F',s')$) l'angle di\'edral (euclidien) du simplexe $s\in \mathcal{T}_{\widetilde M}$ en la face $F$ (resp. du simplexe $s'\in \mathcal{T}_{\widetilde X}$ en la face $F'$). Les nombres $\theta$ et $\theta '$ sont choisis positifs. L'application $\tilde{f_0}$ d'un simplexe $s$ sur un
simplexe $s'$ peut pr\'eserver ou renverser l'orientation (on rappelle
que cette application est affine en restriction \`a $s$) et nous
poserons $\epsilon (s)=\epsilon (s')=\pm 1$ suivant le cas
consid\'er\'e.
 
\begin{lemme}
Soit $F'$ une face de codimension $2$ image de $F$, le degr\'e transverse de $\tilde{f_0}$ en $F$, not\'e $\textrm{deg}_{F}\tilde{f_0}$, v\'erifie,
$$2\pi \textrm{deg}_{F}\tilde{f_0}=\pm \sum_{s'/F'\subset s'}\epsilon (s')\theta '(F',s')\,.$$
\end{lemme}

\begin{preuve}

 Pour $F$ telle que $\tilde{f_0}(F)=F'$ et
$s\in\mathcal{T}_{\widetilde M}$ tels que $F\subset s$,
$\tilde{f_0}(\mathcal{C}\cap s)$ se projette sur $\mathcal{C}'\cap s'$
(o\`u $s'=\tilde{f_0}(s)$) qui est un arc d'angle de valeur absolue
$\theta '(F',s')$. On peut choisir les orientations de
$\widetilde X$ et $F'$ sont telles que l'angle orient\'e de la
projection de $\tilde{f_0}(\mathcal{C}\cap s)$ est $+ \theta '(F',s')$
si $\epsilon (S')=+1$, et $- \theta '(F',s')$ si $\tilde{f_0}$
renverse l'orientation de $s$.  La quantit\'e $\summ_{s'/F'\subset
s'}\epsilon (s')\theta '(F',s')$ repr\'esente donc l'angle orient\'e
total de la projection de $\tilde{f_0}(\mathcal{C})$ sur
$\mathcal{C}'$, c'est-\`a-dire $2\pi \textrm{deg}_{F'}\tilde{f_0}$. Si
l'orientation de $\widetilde X$ est renvers\'ee la relation devient
$2\pi \textrm{deg}_{F'}\tilde{f_0}=-\summ_{s'/F'\subset s'}\epsilon
(s') \theta '(F',s')$.
 
 Consid\'erons alors une d\'eformation de
$\rho _0$, soit $\rho_t$, que nous supposerons $C^1$ en $t$. Nous
construisons l'application $\tilde{f_t}$ de la mani\`ere suivante :

\begin{eqnarray*}
\forall i=1,...,N\,,\, & \tilde{f_t}(m_i) & = y_i \\
\forall \gamma\in\Pi_1(M)\,,\, & \tilde{f_t}(\gamma m_i) & = \rho_t(\gamma)y_i\\
\end{eqnarray*}
et ensuite on \'etend $\tilde{f_t}$ de mani\`ere affine dans chaque
simplexe.  La collection des simplexes images et leurs sommets varient
de mani\`ere $C^1$ en $t$. Nous noterons cette collection
$\mathcal{T}_{\widetilde X}(t)$.  Tous les simplexes de
$\mathcal{T}_{\widetilde X}(t)$ sont non d\'eg\'en\'er\'es, pour $t$
assez petit; en effet, il suffit de n'en consid\'erer qu'un nombre
fini, les autres s'en d\'eduisant par \'equivariance. Notons
\'egalement que, par construction, $\tilde{f_t}$ d\'epend de mani\`ere
$C^1$ en $t$, en particulier, le volume hyperbolique d'un simplexe de
$\mathcal{T}_{\widetilde X}(t)$ est une fonction $C^1$ de $t$. Soit
$F$ une face de codimension $2$ de $\mathcal{T}_{\widetilde M}$ et
$F'(t)$ son image par $\tilde{f_t}$. Pour $s\in
\mathcal{T}_{\widetilde X}(t)$, nous noterons $\theta '(t;F',s')$
l'angle (positif) di\'edral de $s'$ en $F'$. Nous ne mentionnerons pas
la d\'ependance en $t$ des simplexes de $\mathcal{T}_{\widetilde
X}(t)$ et de leurs faces de codimension $2$ s'il n'y a pas
d'ambigu\"{\i}t\'e. Par ailleurs si $s'(t)=\tilde{f_t}(s)$ et $t$ est
assez petit $\epsilon (s'(t))$ ne d\'epend pas de $t$.
\end{preuve}
\begin{lemme}
$$\frac{d}{dt}\big (\sum_{s'/F'\subset s'}\epsilon (s')\theta '(t;F',s')\big )=0$$
\end{lemme}
\begin{preuve}
Pour $t$ assez petit, la face $F'(t)$ est hom\'eomorphe \`a $F'(0)$;
de m\^eme $s'(t)$ est hom\'eomorphe \`a $s'(0)$ si $s'(t)\in
\mathcal{T}_{\wt X}(t)$ et ils sont tous non-d\'eg\'en\'er\'es. Les
voisinages tubulaires de $F'(t)$ et $F'(0)$ sont aussi hom\'eomorphes
et on peut d\'efinir le degr\'e transverse de $\tilde {f_t}$ gr\^ace
\`a $F'(0)$. Alors, par constance du degr\'e par d\'eformation, pour
$t$ assez petit, on a $\deg_{F'(t)}(\tilde {f_t})=\deg_{F'(0)}(\tilde
{f_0})$.
\end{preuve}

Rappelons la formule de Schl\"afli (cf. \cite{Sch}). Soit $s'$ un
simplexe hyperbolique g\'eod\'esique et $F'$ une de ses faces de
codimension $2$; si $s'(t)$ est une d\'eformation de classe $C^1$ de
$s'=s'(0)$, alors
$$\frac{d}{dt}_{|t=0}\vol (s'(t))=-\summ_{F'\subset s'}\frac{d}{dt}_{|t=0}
(\theta '(t;F'(t),s'(t))\vol_{n-2}(F'(t))$$
o\`u $\vol_{n-2}$ d\'esigne le volume $(n-2)$-dimensionnel de la face 
consid\'er\'e.

Pour $\bs \in \mathcal{T}_M$ choisissons un relev\'e $s\in
\mathcal{T}_{\wt M}$; alors $\tilde{f_t}$ identifie de mani\`ere $\ci$
jusqu'au bord $s$ avec un simplexe hyperbolique de $\mathcal{T}_{\wt
X}$. L'\'equivariance de $\tilde{f_t}$ permet de d\'efinir de
mani\`ere unique une m\'etrique hyperbolique sur $\bs$ dont la
collection produit une m\'etrique $\bg (t)$ sur $M$ qui est continue
et hyperbolique par morceaux. En particulier le volume des faces de
codimension $2$ et les angles di\'edraux en celles-ci sont ceux du
simplexe hyperbolique $\tilde{f_t}(s)$. Soit $\omega$ la forme volume
hyperbolique de $\wt X$, alors
\begin{align*}
\vol (\rho _t) & = \int_M\tilde{f_t}^*(\omega ) = 
\summ_{\bs \in \mathcal{T}_M}\int_{\bs} \tilde{f_t}^*(\omega )\\
 & = \summ_{\bs \in \mathcal{T}_M}\epsilon (s)\vol (\tilde{f_t}(s))=
\summ_{\bs \in \mathcal{T}_M} \epsilon (\bs )\vol (\bs ,\bg (t))
\end{align*}
en d\'efinissant $\epsilon (\bs )=\epsilon (s) = \epsilon
(\tilde{f_t}(s))$. Ici on a identifi\'e, par abus de langage,
$\tilde{f_t}^*(\omega )$ avec une forme diff\'erentielle sur $M$
gr\^ace \`a l'\'equivariance de $\tilde{f_t}$. La formule de
Schl\"afli donne,
\begin{align*}
\frac{d}{dt}\vol (\rho_t) &=
\summ_{\bs\in\mathcal{T}_M}\frac{d}{dt}(\epsilon (\bs )\vol(\bs ,\bg
(t)))\\
 &= \summ_{\bs\mathcal{T}_M}\summ_{\bF\subset \bs}\epsilon
(\bs )\frac{d}{dt}(\btheta (t;\bs ,\bF ))\vol _{n-2}(\bF ,\bg (t))
\end{align*}
o\`u $\btheta (t;\bs ,\bF )$ d\'esigne l'angle di\'edral en $\bF$ du
simplexe $\bs$ mesur\'e \`a l'aide de la m\'etrique $\bg (t)$. Il est
\'egal \`a $\theta '(t;\tilde{f_t}(s) ,\tilde{f_t}(F))$ o\`u $s$ et
$F$ sont des relev\'es respectifs de $\bs$ et $\bF$.
$$\frac{d}{dt}\vol (\rho_t)=\summ_{\bF}\Big (\summ_{\bs /\bF\subset
\bs}\epsilon (\bs )\frac{d}{dt}(\btheta (t;\bF ,\bs ))\Big
)\vol_{n-2}(\bF ,\bg (t))$$
 
La quantit\'e entre parenth\`ese peut se calculer sur $M$ ou bien sur
$\wt M$ car elle ne concerne que l'\'etoile d'une face $\bF$; elle
peut \'egalement se calculer sur $\wt X$ par d\'efinition de $\bg
(t)$. Le lemme pr\'ec\'edent montre que, pour toute face $\bF$,
$$\summ_{\bs\supset \bF}\epsilon (\bs )\frac{d}{dt}(\btheta (t; \bF , \bs ))=0$$
Ce qui prouve que $\frac{d}{dt}\vol (\rho_t)=0$.\hfill\finpr

Une cons\'equence imm\'ediate du th\'eor\`eme \ref{vol-const} est le
corollaire suivant. Notons $\mathcal{R}(\Pi_1(M), \Isom (\wt X ))$
l'espace des repr\'esentations du groupe fondamental d'une vari\'et\'e
$M$ dans le groupe d'isom\'etries de l'espace hyperbolique.

\begin{corollaire}
\label{nb-fini}
Soit $M$ une vari\'et\'e diff\'erentielle ferm\'ee et orient\'ee, alors la fonctionnelle,
$$\vol : \mathcal{R}(\Pi_1(M), \Isom (\wt X ))\to \textbf{R}^+$$
prend un nombre fini de valeurs. 
\end{corollaire}
\begin{preuve}
Le groupe d'isom\'etries $\Isom (\wt X)=\PO (n,1)$ est un groupe
alg\'ebrique; par ailleurs, $\Pi_1(M)$ est de pr\'esentation finie
donc $\mathcal{R}(\Pi_1(M), \Isom (\wt X ))$ est une vari\'et\'e
alg\'ebrique (avec singularit\'es) et poss\`ede un nombre fini de
composantes connexes. Le th\'eor\`eme \ref{vol-const} affirme que la
fonctionnelle $\vol$ est constante sur chaque composante connexe.
\end{preuve}
\begin{remarque}
Ce r\'esultat est \'enonc\'e dans \cite{Rez2}, toutefois la preuve est
incompl\`ete sauf, peut-\^etre, en dimension $3$. Celle pr\'esent\'ee
ci-dessus nous a \'et\'e sugg\'er\'ee par J.-P.~Otal (voir \cite{Bon2}).
\end{remarque}
Consid\'erons alors les vari\'et\'es hyperboliques ferm\'ees de
dimension $n$. Un th\'eor\`eme de Wang \cite{Wan}
 affirme que, pour
$n\geq 4$ et $V>0$ le nombre de vari\'et\'es hyperboliques ferm\'ees
de volume inf\'erieur \`a $V$ est fini. Ce r\'esultat est notoirement
faux en dimension $3$ et en dimension $2$. Si $X$ d\'esigne une
vari\'et\'e hyperbolique ferm\'ee et $M$ une vari\'et\'e
diff\'erentielle ferm\'ee, nous dirons (voir \cite{Gro}) que $M$
domine $X$ s'il existe une application continue de degr\'e non nul de
$M$ sur $X$. Le th\'eor\`eme \ref{vol-const} permet de donner une
preuve tr\`es simple du r\'esultat suivant :
\begin{theoreme}[T.~Soma \cite{Som}]
Soit $M$ une vari\'et\'e diff\'erentielle ferm\'ee de dimension $3$,
alors il n'existe qu'un nombre fini de vari\'et\'es hyperboliques de
dimension $3$ ferm\'ees domin\'ees par $M$.
\end{theoreme}
\begin{preuve}
D\'esignons par $f : M\to X$ l'application continue de degr\'e non nul
de $M$ sur $X$, o\`u $X$ est une vari\'et\'e hyperbolique
ferm\'ee. L'application $f$ induit un morphisme $f_* :
\Pi_1(M)\to\Pi_1(X)$, c'est-\`a-dire une repr\'esentation $\rho$ de
$\Pi_1(M)$ dans $\PO (n,1)$. Par d\'efinition du degr\'e de $f$ nous
avons,
$$\vol (\rho )=\deg (f)\vol (X)\,.$$
Par ailleurs, si on munit $M$ d'une m\'etrique riemannienne
quelconque, le th\'eor\`eme \ref{theo-princ} montre que ce volume est
born\'e par un nombre ne d\'ependant que de $M$ (et de cette
m\'etrique) que nous noterons $C(M)$. Nous avons donc,
$$\deg (f) \vol (X)\leq C(M)$$
c'est-\`a-dire, $\deg (f)\leq C(M)/\vol (X)$. Le volume d'une
vari\'et\'e hyperbolique compacte est born\'ee inf\'erieurement par
une constante universelle $v_n$ ne d\'ependant que de la dimension $n$
gr\^ace au lemme de Margulis (voir \cite{Bur-Zal}). En
cons\'equence,
$$\deg (f)\leq C(M)/v_3\,.$$
Il n'y a donc qu'un nombre fini de valeurs possibles pour le degr\'e
de l'application $f$. De m\^eme $\deg (f)\vol (X)=\vol (\rho)$ ne
prend qu'un nombre fini de valeurs d'apr\`es le corollaire
\ref{nb-fini} . Le volume des vari\'et\'es $X$ ferm\'ees domin\'ees
par une vari\'et\'e ferm\'ee fixe $M$ ne peut donc prendre qu'un
nombre fini de valeurs ce qui, d'apr\`es une r\'esultat de W.~Thurston
(\cite{Thu}), montre qu'il ne peut y avoir qu'un nombre fini de telles
vari\'et\'es.
\end{preuve}
Nous terminons en donnant un exemple de telles d\'eformations,
montrant la pertinence du th\'eor\`eme \ref{vol-const}. Il nous a
\'et\'e communiqu\'e par Daryl Cooper par l'interm\'ediaire de Michel
Boileau.
\begin{exemple}(D.~Cooper)
Soit $N$ une vari\'et\'e hyperbolique ferm\'ee de dimension
$3$. Consid\'erons la somme connexe de $N$ avec $S^1\times S^2$,
not\'ee $N\sharp(S^1\times S^2)$, le groupe fondamental de cette
vari\'et\'e est le produit libre $\Pi_1(M)*\textbf{Z}$. Soit $k$ un
noeud homotopiquement nul dans $N\sharp (S^1\times S^2)$ qui rencontre
$S^1\times S^2$ en au moins deux points. D'apr\`es R.~Myers
(\cite{Mye}) on peut trouver de tels noeuds en sorte qu'une chirurgie
de Dehn autour de $k$ transforme $N\sharp (S^1\times S^2)$ en une
vari\'et\'e hyperbolique ferm\'ee $M$ (voir aussi \cite{Boi-Wan} page
797). La proposition 3.2 de \cite{Boi-Wan} permet de construire une
application continue $f : M\to N\sharp (S^1\times S^2)$ de degr\'e
$1$. Par ailleurs il existe \'egalement une application continue,
 $h
: N\sharp (S^1\times S^2)\to N$ de degr\'e $1$ qui consiste \`a
\'ecraser $S^1\times S^2$ en un point. Nous obtenons donc une
application continue de degr\'e $1$,
$$h\circ f:M\longrightarrow N$$
 
et une repr\'esentation
$\rho=h_*\circ f_* : \Pi_1(M)\to Pi_1(N)\subset \PO(3,1)$. Le volume
de cette repr\'esentation est,
$$\vol (\rho )=\vol (N)>0$$
 
car $h\circ f$ est degr\'e $1$ et
l'image de $\rho$ est le groupe fondamental de $N$.
 Par ailleurs
$\rho$ se d\'ecompose en,
$$\rho :
\Pi_1(M)\overset{f_*}{\longrightarrow}\Pi_1(M)*\textbf{Z}\overset{h_*}
{\longrightarrow}\PO(3,1)\,.$$
Le facteur libre $\textbf{Z}$ permet alors de d\'eformer $h_*$ sans
contrainte et donc de produire des d\'eformations non
 triviales (ce
fait est \'el\'ementaire et sa v\'erification est laiss\'ee au
lecteur).
En augmentant le nombre de facteurs $S^1\times S^2$
nous pouvons ais\'ement augmenter le nombre de param\`etres
disponibles pour d\'eformer $\rho$. 
\end{exemple}
\begin{remarque}
Il serait int\'eressant de construire de telles d\'eformations en dimension sup\'erieure ou \'egale \`a $4$. Il est facile d'en construire de volume nul, mais des exemples de volume non nul restent \`a d\'ecrire.
\end{remarque}

\bibliographystyle{plain}
\bibliography{biblio-rep}

\begin{flushleft}
G{\'e}rard BESSON\\
INSTITUT FOURIER\\
Laboratoire de Math\'ematiques\\
UMR5582 (UJF-CNRS)\\
BP 74\\
38402 St MARTIN D'H\`ERES Cedex (France)\\
\smallskip
G.Besson@fourier.ujf-grenoble.fr\\
\medskip
Gilles COURTOIS\\
{\'E}COLE POLYTECHNIQUE\\
Centre de math\'ematiques\\
UMR7640 (CNRS)\\
91128 PALAISEAU Cedex (France)\\
\smallskip
Courtois@math.polytechnique.fr\\
\medskip  
Sylvestre GALLOT\\
INSTITUT FOURIER\\
Laboratoire de Math\'ematiques\\
UMR5582 (UJF-CNRS)\\
BP 74\\
38402 St MARTIN D'H\`ERES Cedex (France)\\
\smallskip
Sylvestre.Gallot@fourier.ujf-grenoble.fr
\end{flushleft}

\end{document}